\documentclass[a4paper,10pt,reqno]{amsart}

\usepackage[T1]{fontenc}
\usepackage[utf8]{inputenc}
\usepackage{lmodern}
\usepackage[english]{babel}

\usepackage{tikz}

\usepackage{xpatch}
\makeatletter
\xpatchcmd{\@settitle}{\uppercasenonmath\@title}{\huge\boldmath}{}{\PatchFailure}
\makeatother

\usepackage[%
	left=2.3cm,       
	right=2.3cm,      
	top=3.3cm,        
	bottom=3.3cm,     
	heightrounded,    
	bindingoffset=0mm 
]{geometry}

\usepackage{amsmath}
\usepackage{amssymb}
\usepackage{mathtools}
\usepackage{mathrsfs}
\usepackage{xfrac}
\numberwithin{equation}{section}

\usepackage[hidelinks]{hyperref} 
\usepackage[noabbrev,capitalize]{cleveref}

\usepackage{lipsum}

\usepackage{bookmark}

\theoremstyle{plain}
	\newtheorem{theorem}{Theorem}[section]
	
	\newtheorem*{lemma*}{Lemma}

\theoremstyle{definition}

	\newtheorem{exercise}[theorem]{Exercise}
	\newtheorem{remark}[theorem]{Remark}

\usepackage[initials]{amsrefs}

\usepackage{pgfmath}

\usepackage{enumerate}

\usepackage{url}

\newcommand{\N}{\mathbb{N}}
\newcommand{\R}{\mathbb{R}}
\newcommand{\Z}{\mathbb{Z}}
\newcommand{\T}{\mathbb{T}}

\newcommand{\eps}{\varepsilon}

\renewcommand{\phi}{\varphi}
\renewcommand{\rho}{\varrho}
\renewcommand{\theta}{\vartheta}

\mathchardef\ordinarycolon\mathcode`\:
\mathcode`\:=\string"8000
\begingroup \catcode`\:=\active
  \gdef:{\mathrel{\mathop\ordinarycolon}}
\endgroup

\def\Xint#1{\mathchoice
{\XXint\displaystyle\textstyle{#1}}%
{\XXint\textstyle\scriptstyle{#1}}%
{\XXint\scriptstyle\scriptscriptstyle{#1}}%
{\XXint\scriptscriptstyle\scriptscriptstyle{#1}}%
\!\int}
\def\XXint#1#2#3{{\setbox0=\hbox{$#1{#2#3}{\int}$ }
\vcenter{\hbox{$#2#3$ }}\kern-.6\wd0}}

\def\aint{\Xint-}

\allowdisplaybreaks

\begin{document}

\title[Introduction to the theory of mixing for incompressible flows]{Introduction to the theory of mixing \\ for incompressible flows}

\author[Gianluca Crippa]{Gianluca Crippa}
\address[G.~Crippa]{Departement Mathematik und Informatik, Universit\"at Basel, Spiegelgasse 1, CH-4051 Basel,\break Switzerland}
\email{gianluca.crippa@unibas.ch}


%
%
%
\begin{abstract}
In these lecture notes, we provide an introduction to the theory of mixing for incompressible flows from a PDE perspective.
We discuss both the Lagrangian (ODE) and Eulerian (PDE, continuity equation) viewpoints, and introduce suitable notions of mixing scales that quantify the degree to which a scalar field transported by a velocity field becomes mixed.
We then address the problem of establishing universal lower bounds on the time evolution of the mixing scale. This is first done in the smooth setting, using energy estimates and flow-based arguments, and later in the Sobolev setting, relying on quantitative estimates for regular Lagrangian flows.
Finally, we present recent results concerning the sharpness of these lower bounds, their implications for the geometry and regularity of regular Lagrangian flows, and connections with more recent developments in the literature.
\end{abstract}

\maketitle

\tableofcontents

\section{Introduction}

Mixing plays a central role in the mathematical understanding of fluid flows, both from a theoretical and an applied perspective. It is also a key concept in many physical phenomena, ranging from atmospheric and oceanic circulation to combustion, chemical reactions, and industrial processes. Moreover, mixing is a main mechanism in the theory of turbulence within fluid dynamics, where it underlies many fundamental processes of energy transfer and dissipation. From a mathematical standpoint, mixing can be studied through a variety of frameworks, including dynamical systems, optimization, and partial differential equations (PDE). In these lecture notes, we focus on the PDE approach. Specifically, we consider the case of passive advection: a passive scalar (such as a pollutant, a dye, or temperature) is transported by a prescribed, incompressible velocity field. A  familiar example is when one pours cream into coffee and stirs, causing coffee and cream to mix. The main objective is to understand the mechanisms and rates by which such a passive scalar homogenizes toward its spatial average. We refer the reader to the survey paper~\cite{MR4733112}, where the connections between mixing and other phenomena in mathematical fluid dynamics are discussed.

In these notes, we will mostly work on the $d$-dimensional torus, which we denote by $\T^d$, and which we take to have unit side length. Our main focus will be the case $d=2$. Working instead on the whole space~$\R^d$ is quite similar, although some aspects, such as scaling analysis and explicit examples, present some differences.  The situation on bounded domains with boundaries is more delicate and requires extra care, so we will not address it here. As usual, $x \in \T^d$ denotes the spatial variable and $t>0$ the time variable.

We denote by $u = u(t,x) : [0,+\infty[ \times \T^d \to \R^d$ a velocity field. We always assume that the velocity field is uniformly bounded in space and time, and that it is divergence free:
\begin{equation}\label{e:divfree}
{\rm div}\, u(t,\cdot) = 0 
\qquad \text{ for all $t \geq 0$.}
\end{equation}
The divergence-free condition~\eqref{e:divfree} encodes the incompressibility of the flow: under the dynamics, the Lebesgue measure of any set is preserved. Appropriate regularity assumptions on the velocity field will be introduced later, as needed.

There are two complementary viewpoints on the problem of passive advection. In the Lagrangian viewpoint (ODE), we consider the flow map $\Phi=\Phi(t,x) : [0,+\infty[ \times \T^d \to \T^d$ generated by the velocity field, defined by
\begin{equation}\label{e:ODE}\tag{ODE}
\left\{ \begin{array}{l}
\dot \Phi (t,x) = u (t,\Phi(t,x)) \\
\Phi(0,x) = x \,. 
\end{array}\right.
\end{equation}
In this formulation, we follow the trajectory of each individual fluid particle. In the Eulerian viewpoint~(PDE), we consider the advection of a passive scalar $\rho = \rho(t,x) : [0,+\infty[ \times \T^d \to \R$ governed by the continuity equation
\begin{equation}\label{e:PDE}\tag{CE}
\left\{ \begin{array}{l}
\partial_t \rho + {\rm div}\, (u\rho) = 0 \\
\rho(t=0) = \bar{\rho} = \bar{\rho}(x) \,,
\end{array}\right.
\end{equation}
which, under the divergence-free constraint~\eqref{e:divfree} on the velocity field, is equivalent to the transport equation $\partial_t \rho + u \cdot \nabla \rho = 0$. This equation models the evolution of a  profile under the action of the prescribed velocity field. Physically, it expresses the local conservation of the mass of the advected quantity. We also recall that the two viewpoints (Lagrangian and Eulerian) are connected by the theory of characteristics, which provides a formula in terms of the flow:
\begin{equation}\label{e:character}
\rho(t,x) = \bar\rho\big( \Phi^{-1}(t,\cdot)(x) \big) 
\qquad \text{ for all $x \in \T^d$ and $t \geq 0$.}
\end{equation}
\begin{exercise}
Show that formula~\eqref{e:character} provides a solution to~\eqref{e:PDE}. 
\end{exercise}

For the initial condition $\bar\rho$ of the passive scalar, we assume it has zero average, that is
\begin{equation}\label{e:zeroaver}
\aint_{\T^d} \bar\rho(x) \, dx = 0 \,.
\end{equation}
\begin{exercise}
Show that the zero-average condition~\eqref{e:zeroaver} is propagaged in time, that is, for every solution~$\rho$ of~\eqref{e:PDE}, it holds
$$
\aint_{\T^d} \rho(t,x) \, dx = 0 
\qquad \text{ for all $t \geq 0$.}
$$
\end{exercise}
Therefore, the zero-average condition~\eqref{e:zeroaver} is not a real restriction: since constants are preserved under divergence-free advection, we can always reduce to this case by subtracting the average of the initial datum. This can be seen as a gauge condition, which is sometimes technically convenient.


Our objective is to study the mixing of the passive scalar $\rho$ solution of~\eqref{e:PDE}, namely its convergence to zero (i.e., to its spatial average). Ideally, this occurs through stirring and the creation of fine filamentation. 
A natural  idea would be to study the variance of~$\rho$ as a function of time. Thanks to the zero-average condition~\eqref{e:zeroaver}, the variance simply coincides with the square of the $L^2$-norm of $\rho$. However, at least formally (that is, for sufficiently smooth solutions), every $L^p$-norm of $\rho$ is conserved in time. In particular, the variance remains constant. Therefore, the study of the variance does not provide any useful information about the mixing of the passive scalar.
\begin{exercise}
Show that, for sufficiently smooth solutions of~\eqref{e:PDE}, any $L^p$-norm is constant in time. You can use the representation formula~\eqref{e:character}, or manipulate directly the first equation in~\eqref{e:PDE}. 
\end{exercise}
\begin{remark}
The variance (or the $L^2$-norm) is a suitable quantity to study mixing in the presence of diffusion, that is, for solutions to the advection-diffusion equation $\partial_t \rho + {\rm div}\, (u\rho) = \Delta \rho$. In these notes we will not address this case and we refer the interested reader for instance to~\cite{MR4733112} and the references therein. 
\end{remark}
Mixing occurs through filamentation and the intertwining of level sets of the passive scalar: the scalar is stretched into thin, intertwined filaments by the flow, causing regions of different values to blend together more efficiently. We refer for instance to~\cite{MR3824693} for an introduction to the concepts of stirring, mixing, and transport. A suitable way to capture this phenomenon is to consider weak convergence, which expresses the fact that local averages of $\rho$ converge to zero (the spatial average of the passive scalar). Weak convergence to zero is the property that for every smooth test function $\varphi = \varphi(x) \in C^\infty(\T^d)$,
$$
\int_{\T^d} \rho(t,y)\,\varphi(y)\,dy \to 0
\qquad \text{as } t \to \infty\,,
$$
or alternatively in our case, for every ball~$B(x,r)\subset\T^d$,
$$
\int_{B(x,r)} \rho(t,y) \,dy \to 0
\qquad \text{as } t \to \infty\,.
$$

Our objective is to formulate a fully quantitative theory of mixing in physically relevant settings. We do not only want to state that a passive scalar gets mixed, but rather to introduce a precise notion of mixing scale (a characteristic length below which the passive scalar appears approximately homogeneous) and to establish quantitative dynamical bounds on the evolution of this mixing scale in time.

In order to achieve this program, the following three main steps need to be addressed:
\begin{itemize}
\item[(1)] Definition of mixing scales. We begin by introducing in Section~\ref{s:scales} two notions of mixing scale: a geometric one and a functional one. 
\item[(2)] Universal lower bounds. We establish universal lower bounds on the mixing scale, under natural energetic constraints on the velocity fields. We do this via energy estimates (in Section~\ref{s:eneest}), through Lagrangian estimates for Lipschitz velocity field (in Section~\ref{s:flowestimates}), and by exploiting the theory of quantitative estimates for regular Lagrangian flows for Sobolev velocity fields from~\cite{MR2369485} (in Sections~\ref{s:CDLproof}, after some preliminaries in Sections~\ref{s:functplan}~and~\ref{s:hatools}).  
\item[(3)] Sharpness of the bounds. Are there configurations for which these lower bounds are attained? We present a combinatorial mixing scheme based on ``slice-and-dice'' mechanism in Section~\ref{s:bressan}, and describe some more results and give references in Section~\ref{s:saturation}.
\end{itemize}

\begin{remark}
Concerning step~(2), we make the following observations. If no energetic constraint is imposed on the velocity field, no universal lower bound can hold. Indeed, by rescaling the velocity field, one can produce arbitrarily fast mixing. Under energetic constraints, we can establish lower bounds on the mixing scale, but not upper bounds. This is because such constraints allow for flows that do not mix at all (for instance,  translations on the torus). This situation contrasts with the approach in dynamical systems, where one assumes geometric or structural properties of the flow (such as the presence of hyperbolic points), and then proves that every passive scalar is mixed. Technically, in that setting one deals with the decay of correlations for observables. The PDE and dynamical-systems approaches are distinct, but nevertheless deeply connected.
\end{remark}

Along the way, we will also encounter aspects of the DiPerna-Lions-Ambrosio~\cite{MR1022305,MR2096794} theory of flows associated with non-Lipschitz velocity fields (Sobolev or bounded variation), see Section~\ref{s:bressan}. This theory has had profound consequences for the study of nonlinear PDEs in mathematical physics, including kinetic equations, fluid dynamics, and conservation laws. The constructions involved in step~(3) above provide an important perspective on the regularity and irregularity of such flows, and we will return to this point in Section~\ref{s:saturation}. 

\medskip

{\bf Notation.} We mostly employ standard notation throughout these notes. We denote by $B(x,r)$ the ball centered at $x$ with radius $r$. Given a measurable set $A$, we write $|A|$ for its Lebesgue measure. The symbol $\aint_A$ denotes the average integral over $A$. We use the standard notation $L^p$, $W^{s,p}$, and $H^s$ for Lebesgue and Sobolev spaces, and similarly for their time-dependent counterparts. The arrow $\rightharpoonup$ denotes weak convergence. The symbol $\lesssim$ indicates an inequality that holds up to a multiplicative constant depending only on fixed parameters appearing in the argument. In some vector-calculus computations, summation over repeated indices is understood.

\medskip

{\bf Acknowledgments.} These notes have evolved from a series of advanced lectures I had the opportunity to deliver on several occasions. I wish to thank the colleagues who invited me to present these lectures, the institutions that hosted and supported them, and the participants for their engagement and constructive feedback. I also acknowledge support by the ERC Starting Grant 676675 FLIRT, by the SNSF Project 212573 FLUTURA, and by the SPP 2410 CoScaRa funded by the DFG  through the project 217527 funded by the SNSF.

\section{The geometric mixing scale and the functional mixing scale}\label{s:scales}
\begin{figure}[h]
\begin{center}
\includegraphics[scale=0.6]{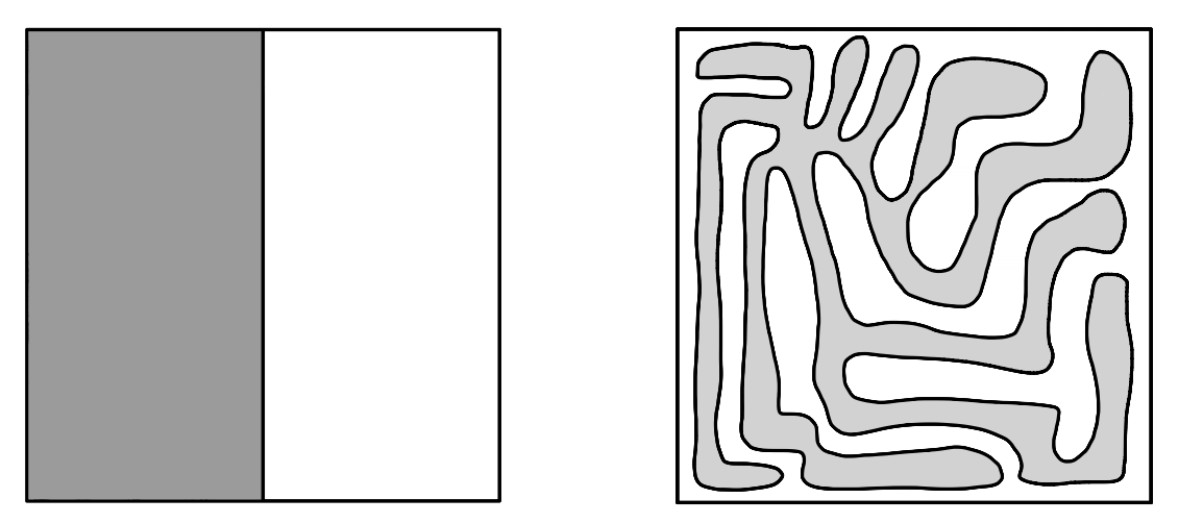}
\caption{An initial configuration, and a more mixed final configuration.}\label{f:mixscale_1}
\end{center}
\end{figure}

If we look at the initial datum (left) and the final datum (right) in Figure~\ref{f:mixscale_1}, it is clear that the final configuration is (much) more mixed than the initial one. However, our goal is to make this observation quantitative: what is the mixing scale of these two configurations?

\begin{figure}[h]
\begin{center}
\includegraphics[scale=0.6]{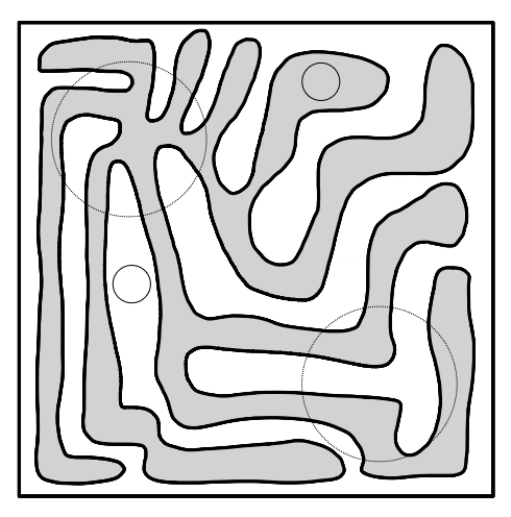}
\caption{Finding the optimal resolution at which the two phases are balanced.}\label{f:mixscale_2}
\end{center}
\end{figure}

The mixing scale should be understood as the characteristic length at which the configuration appears approximately homogeneous. Put differently, if as observers we only resolve the system at a limited resolution, the mixing scale is the finest resolution at which the different phases of the fluid can no longer be distinguished. For instance, in Figure~\ref{f:mixscale_2}, the circles with solid boundary are too small: at that scale the two phases of the fluid remain clearly visible. The circles with dotted boundary have roughly the size of the mixing scale, since within each circle of such radius we see an essential balance between the two phases. Any circle with greater diameter, on the other hand, would be too large: the fluid looks mixed at that scale, but it is already mixed at finer scales as well. We would like to define the radius of the dotted circles as the mixing scale of the passive scalar. 

There are two main ways of formalizing this idea: the geometric mixing scale ${\rm mix}_{\rm g}(\rho)$, and the functional mixing scale ${\rm mix}_{\rm f}(\rho)$. Since the mixing scales depend only on the spatial configuration, we define them first for a time-independent $\rho=\rho(x)$. Later, we will study how this quantity evolves in time under the dynamics. The two notions of mixing scale are related and behave in similar ways, but they are not equivalent (see~\cite{MR3026556,MR4026549}).

\subsection{The geometric mixing scale}\label{ss:defgeomix} The geometric mixing scale (introduced in~\cite{MR2033002}) is based on the idea of the intertwining of the level sets of the passive scalar. Consider the case of a binary configuration $\rho \in \{ \pm 1\}$ and fix an accuracy parameter $0 < \kappa < \sfrac{1}{2}$. The geometric mixing scale ${\rm mix}_{\rm g} (\rho)$ is defined as the infimum of all~$\eps>0$ such that 
\begin{equation}\label{e:mixgeocond}
\kappa \leq \frac{\left| \{ \rho=+1\} \cap B(x,\eps) \right|}{|B(x,\eps)|} \leq 1-\kappa
\qquad \text{ for all $x \in \T^d$.}
\end{equation}
In words, for example in the case $\kappa = \sfrac{1}{3}$, condition~\eqref{e:mixgeocond} requires that in every ball of radius~$\eps$, at least one third of the volume is occupied by the $+1$ phase and at least one third by the $-1$ phase. This definition can be generalized to non-binary configurations, by replacing~\eqref{e:mixgeocond} with
$$
\frac{1}{\| \rho\|_{L^\infty(\T^d)}} \left| \, \aint_{B(x,\eps)} \rho \, dy \right| \leq \kappa' 
\qquad \text{ for all $x \in \T^d$,}
$$
where $0<\kappa'<1$. 

\begin{exercise}
Find the relation between $\kappa'$ and $\kappa$ in the case $\rho \in \{\pm 1\}$. 
\end{exercise}

\begin{exercise}
Show that, if ${\rm mix}_{\rm g} (\rho)=0$, then $\rho=0$ almost everywhere. 
\end{exercise}

\begin{exercise}
Let $\{\rho_n\}$ be equi-bounded in $L^\infty$ and with zero average. Show that:
\begin{itemize}
\item[(i)] If ${\rm mix}_{\rm g} (\rho_n)\to0$ for a given $\kappa'$, then not necessarily $\rho_n \rightharpoonup 0$ weakly in $L^2$. 
\item[(ii)] ${\rm mix}_{\rm g} (\rho_n )\to 0$ for all $\kappa'>0$ if and only if $\rho_n \rightharpoonup 0$ weakly in $L^2$. 
\end{itemize}
\end{exercise}

\subsection{The functional mixing scale} The concept of functional mixing scale (introduced in this context in~\cite{MR2801050,MR2200439}, see also the review paper~\cite{MR2876867}) is motivated by the idea that mixing corresponds to the creation of small scales, which in turn correspond to high frequencies. Denoting by~$\hat{\rho}(k)$, for $k \in \Z^d$, the Fourier coefficients of $\rho$, we define
\begin{equation}\label{e:fctmix}
{\rm mix}_{\rm f} (\rho) = \left[\sum_{\substack{k \in \Z^d \\ k\not = 0}} \frac{1}{|k|^2} \, | \hat{\rho}(k)|^2\right]^{\sfrac{1}{2}} \,.
\end{equation}
\begin{exercise}
Show that
$$
\hat{\rho}(0) = c \int_{\T^d} \rho(y) \, dy \,,
$$
where the constant $c$ depends on the chosen normalization in the definition of the Fourier coefficients. Observe as a consequence that the zero-average condition~\eqref{e:zeroaver} guarantees that in~\eqref{e:fctmix} it does not really matter that we did not include the term $k=0$ in the sum.
\end{exercise}

Note that Plancherel's theorem states that 
$$
\| \rho \|_{L^2(\T^d)} = \left[\sum_{\substack{k \in \Z^d \\ k\not = 0}}  | \hat{\rho}(k)|^2\right]^{\sfrac{1}{2}} \,.
$$
Recall that the $L^2$-norm of a solution $\rho(t,x)$ to~\eqref{e:PDE} is conserved in time. Therefore, the fact that ${\rm mix}_{\rm f} (\rho(t,\cdot)) \to 0$ for a solution expresses that energy is being transferred from low to high frequencies. Indeed, for large values of $|k|$, the corresponding term in the sum in~\eqref{e:fctmix} contributes less to the total value, so the decay of the functional mixing scale reflects the progressive shift of energy toward higher frequencies.

\begin{exercise}\label{ex:funct}
Let $\{\rho_n\}$ be equi-bounded in $L^2(\T^d)$ and with zero average. Show that:
$$
{\rm mix}_{\rm f} (\rho_n) \to 0  
\qquad \Longleftrightarrow \qquad
\rho_n \rightharpoonup 0 \quad \text{weakly in $L^2(\T^d)$}. 
$$
\end{exercise}

\subsection{A brief overview of fractional Sobolev spaces} Recalling that the Fourier transform converts differentiation in physical space into multiplication by the Fourier variable, and using Plancherel’s theorem, we can interpret the right-hand side of~\eqref{e:fctmix} as the $L^2$-norm of a derivative of order $-1$ of $\rho$, that is, as the (homogeneous) Sobolev norm $\|\rho\|_{\dot{H}^{-1}}$ of negative order $-1$. More in general, for any $s \in \R$ we can define
\begin{equation}\label{e:defsob}
\| f \|_{\dot{H}^s(\T^d)}^2 = \sum_{\substack{k \in \Z^d \\ k\not = 0}} |k|^{2s} \, | \hat{f}(k)|^2 \,.
\end{equation}

The quantity in~\eqref{e:defsob} can be used to define homogeneous fractional Sobolev spaces. The term homogeneous (indicated by the dot in the notation) means that the definition involves only the integral norm of the derivative of the given order, and not the integral norm of the function itself. In general, the rigorous definition of these spaces is somewhat delicate: they must be introduced as suitable quotients of spaces of distributions. However, we will not enter into the full functional-analytic theory. Instead, we will use these norms simply as a way to measure the size of ``nice enough'' functions.

We now motivate the choice of the exponent $-1$ in the definition of the functional mixing scale in~\eqref{e:fctmix}. In fact, the statement in Exercise~\ref{ex:funct} would remain valid for any $s<0$. However, the $\dot{H}^{-1}$-norm on the torus defines a characteristic length, which can be interpreted as the length scale at which the configuration is mixed, analogous to the length given by the geometric mixing scale. We will now explain this by means of a scaling computation. 

\begin{exercise}
Understand what could be a definition in Fourier variable of the homogeneous Sobolev space $\dot{W}^{s,p}(\T^d)$, for $s \in \R$ and $1<p<\infty$. In fact, we will not really need the precise rigorous definition, but rather the operative properties of these spaces, so we will not focus too much on formal definitions.
\end{exercise}

If $f=f(x)$ is defined on $\T^d$ and $\sfrac{1}{\lambda}\in\N$, then $f_\lambda(x) = f(\sfrac{x}{\lambda})$ is also $1$-periodic, and hence well-defined as a function on $\T^d$. Denoting by $\nabla^s$ the derivative of order $s$, we formally have the relation 
$$
\nabla^s \left( f\left(\frac{x}{\lambda}\right)\right) = \frac{1}{\lambda^s} \left(\nabla^s f\right) \left(\frac{x}{\lambda}\right) \,.
$$
Hence we can formally compute
$$
\begin{aligned}
\| f_\lambda \|_{\dot{W}^{s,p}(\T^d)}^p
&= \int_{\T^d} | \nabla^s f_\lambda|^p (y) \, dy
= \int_{\T^d} \frac{1}{\lambda^{sp}} | \nabla^s f|^p \left( \frac{y}{\lambda}\right) \, dy \\
&= \frac{1}{\lambda^{sp}} \, \frac{1}{\lambda^d} \int_{\T^d} | \nabla^s f|^p (x) \, \lambda^d \, dx
= \lambda^{-sp} \| \nabla^s f \|_{L^p(\T^d)}^p \,,
\end{aligned}
$$
(notice that the periodicity of the domain of integration gives the factor $\sfrac{1}{\lambda^d}$ in the change of variable). Hence, we obtain
\begin{equation}\label{e:resctorus}
\| f_\lambda \|_{\dot{W}^{s,p}(\T^d)} = \lambda^{-s} \| f \|_{\dot{W}^{s,p}(\T^d)} \,.
\end{equation}
In particular, the $\dot{H}^{-1}$-norm scales as a length on the torus and therefore it determines a characteristic length scale of the given configuration.

\begin{exercise} Show formally that on $\R^d$ the scaling relation
\begin{equation}\label{e:rescspace}
\| f_\lambda \|_{\dot{W}^{s,p}(\R^d)} = \lambda^{\frac{d}{p}-s} \| f \|_{\dot{W}^{s,p}(\R^d)}
\end{equation}
holds. On $\R^d$, there is no factor $\sfrac{1}{\lambda^d}$ appearing in the change of variables. We can heuristically explain the difference between $\R^d$ and $\T^d$ as follows: when rescaling a single bump on $\R^d$, we still obtain a single bump, whereas on $\T^d$ the rescaling produces multiple copies of the bump, due to periodicity.
\end{exercise}

\begin{exercise}
Prove rigorously~\eqref{e:resctorus} and~\eqref{e:rescspace} by relying on properties of the Fourier transform.  
\end{exercise}

\begin{exercise}
Show that the geometric mixing scale also scales as a length. 
\end{exercise}

\section{Energy estimates for the functional mixing scale}\label{s:eneest}

In this section we begin deriving lower bounds for the mixing scale. We follow~\cite{MR2801050} (see also~\cite{MR3026556}), where such bounds are obtained through energy estimates on the Eulerian side, namely for~\eqref{e:PDE}. These estimates are tailored to the functional mixing scale and have the advantage of being relatively straightforward and direct to derive. However, in certain cases they are not optimal. In what follows, we will consider three different types of energy budgets, one in each of the next subsections.

\subsection{Velocity fields with uniformly-in-time Lipschitz regularity}\label{ss:lipPDE}
Let us assume that the velocity field $u$ has uniform-in-time Lipschitz regularity, that is
\begin{equation}\label{e:vellip}
{\rm Lip}\, (u(t,\cdot)) \leq L \qquad \forall \, t \geq 0 
\end{equation}
for some constant $L$, where we denote by ${\rm Lip}\, (\cdot)$ the Lipschitz constant of a given function or vector field.  
Thanks to the zero-average condition~\eqref{e:zeroaver}, we can consider the unique $\varphi$ with zero average such that~$\Delta \varphi = \rho$ and multiply~\eqref{e:PDE} by $\varphi$. We get
$$
\varphi \, \partial_t \Delta \varphi + \varphi \, {\rm div} \, (u \Delta \varphi) = 0 \,.
$$
We now integrate over $\T^d$ and manipulate the two terms as follows. For the first term we obtain
$$
\begin{aligned}
\int_{\T^d} \varphi \, \partial_t \Delta \varphi \, dx  &= \int_{\T^d} \varphi \, \Delta \partial_t \varphi \, dx 
= - \int_{\T^d} \nabla \varphi \cdot \nabla \partial_t \varphi \, dx  \\
&= -\int_{\T^d} \nabla \varphi \cdot \partial_t \nabla \varphi \, dx = - \frac{1}{2} \frac{d}{dt} \int_{\T^d} | \nabla \varphi |^2 \, dx \,,
\end{aligned}
$$
while the second term can be rewritten as
$$
\begin{aligned}
\int_{\T^d} \varphi \, {\rm div}\, (u \Delta \varphi) \,  dx &= \int_{\T^d} \varphi \, \partial_i (u^i \partial_{jj} \varphi) \, dx 
= - \int_{\T^d} \partial_i \varphi \, u^i \, \partial_{jj} \varphi \, dx \\
&= \int_{\T^d} \partial_j u^i \, \partial_i \varphi \, \partial_j \varphi \, dx
+ \int_{\T^d} u^i \, \partial_i \partial_j \varphi \, \partial_j \varphi \, dx \\
&= \int_{\T^d}  Du :  ( \nabla \varphi \otimes \nabla \varphi) \, dx
+ \frac{1}{2} \int_{\T^d} u^i \, \partial_i ( | \partial_j \varphi|^2 ) \, dx \\
&= \int_{\T^d}  Du :  ( \nabla \varphi \otimes \nabla \varphi) \, dx
- \frac{1}{2} \int_{\T^d} \underbrace{\partial_i u^i}_{={\rm div}\, u}  | \partial_j \varphi|^2  \, dx \\
&= \int_{\T^d}  Du :  ( \nabla \varphi \otimes \nabla \varphi) \, dx \,.
\end{aligned}
$$
Therefore, by using~\eqref{e:vellip}, we can bound
$$
\begin{aligned}
\frac{d}{dt} \int_{\T^d} |\nabla \varphi|^2 \, dx &= 2 \int_{\T^d}  Du :  ( \nabla \varphi \otimes \nabla \varphi) \, dx \\
&\geq -2L \int_{\T^d} |\nabla \varphi|^2 \, dx \,. 
\end{aligned}
$$
Observing that 
$$
\int_{\T^d} | \nabla \varphi |^2 \, dx = \| \rho \|_{\dot{H}^{-1}(\T^d)}^2 \,,
$$
if we integrate the above differential inequality we obtain
\begin{equation}\label{e:firstfct}
\| \rho(t,\cdot) \|_{\dot{H}^{-1}(\T^d)}
\geq
\| \bar\rho\|_{\dot{H}^{-1}(\T^d)} \, e^{-tL} \,.
\end{equation}
The functional mixing scale can decay at most at a negative exponential rate in time. The coefficient in the exponential is determined by the Lipschitz constant of the velocity field, while the lower bound depends linearly on the functional mixing scale of the initial datum. In this case, the derived lower bound has the structure
$$
{\rm mix}_{\rm f} ( \rho(t,\cdot) ) \geq F (t \, , \, \text{energy of $u$}) \cdot G (\bar\rho) \,,
$$
where $F$ is a universal function, independent of the particular solution, while $G$ depends on the initial datum only through its functional mixing scale.

Let us remark that assumptions on the Lipschitz regularity of the velocity field are usually too strong for applications in fluid dynamics; therefore, in the next two subsections we will weaken this assumption.

\subsection{Velocity fields with uniformly-in-time bounded kinetic energy}\label{ss:kinPDE} Let us now assume that the velocity field $u$ has uniform-in-time bounded kinetic energy, that is
\begin{equation}\label{e:velkine}
\| u(t,\cdot) \|_{L^2(\T^d)} \leq K \qquad \forall \, t \geq 0 
\end{equation}
for some constant $K$. We proceed as in~Subsection~\ref{ss:lipPDE}, but we exploit a different rewriting for the second term:
$$
\int_{\T^d} \varphi \, {\rm div}\, (u \Delta \varphi) \,  dx
= 
- \int_{\T^d} u \, \Delta \varphi \cdot \nabla \varphi \, dx \,.
$$
Therefore, using~\eqref{e:velkine} we obtain 
$$
\begin{aligned}
\frac{d}{dt} \int_{\T^d} |\nabla \varphi|^2 \, dx &= -2  \int_{\T^d} u \, \Delta \varphi \cdot \nabla \varphi \, dx 
\geq -2 \| u \|_{L^2(\T^d)} \| \Delta \varphi \|_{L^\infty(\T^d)} \| \nabla \varphi \|_{L^2(\T^d)} \\
&\geq -2 K \| \rho \|_{L^\infty(\T^d)} \| \nabla \varphi \|_{L^2(\T^d)}  \,,
\end{aligned}
$$
which in turn gives
$$
\frac{d}{dt} \| \rho (t,\cdot) \|_{\dot{H}^{-1}(\T^d)}
\geq
-K \| \rho(t,\cdot) \|_{L^\infty(\T^d)}
= -K \| \bar\rho \|_{L^\infty(\T^d)} \,.
$$
As a consequence, we obtain a linear-in-time lower bound on the functional mixing scale
$$
\| \rho (t,\cdot) \|_{\dot{H}^{-1}(\T^d)}
\geq 
\| \bar\rho \|_{\dot{H}^{-1}(\T^d)}
- \left( K \| \bar\rho \|_{L^\infty(\T^d)} \right) t \,.
$$
This, in principle, allows for $\| \rho (t,\cdot) \|_{\dot{H}^{-1}(\T^d)}=0$ and therefore $\rho (t,\cdot)=0$ at some finite time. Such a phenomenon is called perfect mixing in finite time: the solution becomes completely mixed at that finite time. By time-reversal, perfect mixing in finite time also implies nonuniqueness for the initial-value problem for~\eqref{e:PDE}. In Section~\ref{s:bressan} we will see, through an example, that under a uniform-in-time bound on the kinetic energy as~\eqref{e:velkine} both perfect mixing in finite time and nonuniqueness can indeed occur, and we will comment on the connections with the DiPerna-Lions-Ambrosio~\cite{MR1022305,MR2096794} theory of flows of non-Lipschitz velocity fields.

\subsection{Velocity fields with uniformly-in-time bounded enstrophy}\label{ss:ensPDE} Finally, we consider the case in which the velocity field $u$ has uniform-in-time bounded enstrophy, that is
\begin{equation}\label{e:velenstr}
\|u(t,\cdot)\|_{\dot{H}^1} \leq E \qquad \forall \, t \geq 0 
\end{equation}
for some constant $E$. In this case, we rewrite the second term as
\begin{equation}\label{e:earlier}
\begin{aligned}
\int_{\T^d} \varphi \, {\rm div}\, (u \Delta \varphi) \,  dx
&= 
\int_{\T^d} Du : (\nabla \varphi \otimes \nabla \varphi) \, dx
=
\int_{\T^d} \partial_j u^i \, \partial_i \varphi \, \partial_j \varphi \, dx \\
&=
- \int_{\T^d} \partial_j \underbrace{\partial_i u^i}_{={\rm div}\, u} \, \varphi \, \partial_j \varphi \, dx
- \int_{\T^d} \partial_j u^i \, \varphi \, \partial_{ij} \varphi \, dx \\
&=
 - \int_{\T^d} \partial_j u^i \, \varphi \, \partial_{ij} \varphi \, dx \,.
\end{aligned}
\end{equation}
Recalling that $\rho = \Delta \varphi$ and observing that
$$
\| \partial_{ij} \varphi \|_{L^2(\T^d)} \lesssim \| \Delta \varphi \|_{L^2(\T^d)} = \| \rho \|_{L^2(\T^d)} \,,
$$
using~\eqref{e:velenstr} we obtain
$$
\begin{aligned}
\frac{d}{dt} \int_{\T^d} |\nabla \varphi|^2 \, dx &= - 2 \int_{\T^d} \partial_j u^i \, \varphi \, \partial_{ij} \varphi \, dx 
\gtrsim  -2 E  \| \varphi  \|_{L^\infty} \| \rho  \|_{L^2(\T^d)} \\
&=  -2 E  \| \varphi  \|_{L^\infty(\T^d)} \| \bar\rho \|_{L^2(\T^d)} \,.
\end{aligned}
$$
We need to close this estimate, but the $L^\infty$-norm of $\varphi$ is not related to any norm of $\rho$ which appears in the estimates. If we could control
\begin{equation}\label{e:wrong}
\| \varphi \|_{L^\infty(\T^d)} \lesssim \| \nabla \varphi \|_{L^2(\T^d)} 
\qquad\text{(this estimate does not hold!)}
\end{equation}
we would immediately get in any spatial dimension
\begin{equation}\label{e:subopt}
\frac{d}{dt} \| \rho (t,\cdot) \|_{\dot{H}^{-1}(\T^d)}
\gtrsim
-E \| \bar\rho \|_{L^2(\T^d)} \,.
\end{equation}
However, in dimension $d \geq 2$, estimate~\eqref{e:wrong} does not hold. 

\begin{exercise}
Derive  lower bounds on the functional mixing scale under the assumption~\eqref{e:velenstr}. In two dimensions, this can be done by going back to~\eqref{e:earlier} and using Gagliardo-Nirenberg's inequality
$$
\| f \|_{L^4(\T^2)} \lesssim  \| \nabla f \|^{\sfrac{1}{2}}_{L^2(\T^2)} \| f \|^{\sfrac{1}{2}}_{L^2(\T^2)} \,,
$$
while in three dimensions one can use Agmon's inequality 
$$
\| f \|_{L^\infty(\T^3)} \lesssim \| \nabla f \|_{L^2(\T^3)}^{\sfrac{1}{2}} \| \Delta f \|_{L^2(\T^3)}^{\sfrac{1}{2}} \,.
$$
In two dimensions, this gives the lower bound~\eqref{e:subopt}. However,~\eqref{e:subopt} does not hold in higher dimension, at least not as a consequence of these energy estimates. 
\end{exercise} 

The lower bound~\eqref{e:subopt} (valid in two dimensions) 
only yields a linear-in-time lower bound on the functional mixing scale (by the same computations as in Subsection~\ref{ss:kinPDE}), which is not sufficient to rule out perfect mixing in finite time and nonuniqueness for the initial-value problem for~\eqref{e:PDE}. Nevertheless, we will see in Section~\ref{s:CDLproof} that, differently from the situation in Subsection~\ref{ss:kinPDE},  under the assumption~\eqref{e:velenstr} uniqueness for the initial-value problem for~\eqref{e:PDE} is guaranteed as a consequence of the DiPerna-Lions-Ambrosio~\cite{MR1022305,MR2096794} theory, and therefore perfect mixing in finite time cannot occur. This shows that the lower bound~\eqref{e:subopt} obtained through energy estimates should not be expected to be optimal.


\section{Estimates using the flow for uniformly-in-time Lipschitz velocity fields}\label{s:flowestimates}

In this section, we show how to obtain estimates from a Lagrangian point of view, that is, by relying on the formulation~\eqref{e:ODE}. Denoting by~$\Phi$ the flow of the velocity field $u$, then under the assumption ${\rm Lip}\, (u(t,\cdot)) \leq L$ for all $t\geq 0$ we have
\begin{equation}\label{e:regLipflow}
e^{-Lt} \leq \frac{| \Phi(t,x)-\Phi(t,y) |}{|x-y|} \leq e^{Lt} 
\qquad \text{for all $t \geq 0$ and $x,y\in\T^d$}
\end{equation}
as well as
\begin{equation}\label{e:regLipflowinverse}
e^{-Lt} \leq \frac{| \Phi(t,\cdot)^{-1}(x)-\Phi(t,\cdot)^{-1}(y) |}{|x-y|} \leq e^{Lt} 
\qquad \text{for all $t \geq 0$ and $x,y\in\T^d$.}
\end{equation}

\begin{exercise}
Prove~\eqref{e:regLipflow} and~\eqref{e:regLipflowinverse} by observing that
$$
\left| \frac{d}{dt} | \Phi(t,x)-\Phi(t,y)|^2 \right| \leq 2L | \Phi(t,x)-\Phi(t,y)|^2
$$
and then using Gr\"onwall inequality. 
\end{exercise}

\subsection{Lower bound for the functional mixing scale}\label{ss:regestinter}

Recall the representation formula~\eqref{e:character} for the solution of~\eqref{e:PDE} in terms of the flow. If we now assume $\bar\rho \in \dot{H}^1(\T^d)$, we can prove an estimate for the propagation of the $\dot{H}^1$-regularity of the solution as follows:
$$
\begin{aligned}
\| \rho(t,\cdot) \|_{\dot{H}^1(\T^d)}^2 &= \int_{\T^d} | \nabla \rho(t,x)|^2 \, dx = \int_{\T^d} | \nabla (\bar \rho (\Phi^{-1}(t,\cdot)(x)))|^2 \, dx \\
&\leq \int_{\T^d} | (\nabla \bar\rho) (\Phi^{-1}(t,\cdot)(x)) |^2 \, |\nabla (\Phi^{-1}(t,\cdot)(x)) |^2 \, dx \\
&\leq ({\rm Lip}\, \Phi^{-1})^2 \int_{\T^d} | \nabla\bar\rho|^2 \, dx \leq e^{2Lt} \, \| \bar\rho\|_{\dot{H}^1(\T^d)}^2 \,,
\end{aligned}
$$
that is,
\begin{equation}\label{e:propagPDE}
\| \rho(t,\cdot) \|_{\dot{H}^1(\T^d)}
\leq 
\| \bar\rho\|_{\dot{H}^1(\T^d)} \, e^{Lt}  \,.  
\end{equation}
Estimate~\eqref{e:propagPDE} expresses that, for a uniformly-in-time Lipschitz velocity field, if the initial datum belongs to $\dot{H}^1(\T^d)$, then the unique solution remains in $\dot{H}^1(\T^d)$ for all times, with at most exponential-in-time growth of the corresponding norm. 
\begin{exercise}
Prove the interpolation inequality
\begin{equation}\label{e:interpol}
\| f \|_{L^2(\T^d)}^2 \leq \| f \|_{\dot{H}^{-1}(\T^d)} \, \| f \|_{\dot{H}^1(\T^d)} 
\qquad \text{ for all $f \in \dot{H}^1(\T^d)$.}
\end{equation}
\end{exercise}
By using~\eqref{e:propagPDE} and~\eqref{e:interpol} we can conclude that
\begin{equation}\label{e:lowerfunctinter}
{\rm mix}_{\rm f} (\rho(t,\cdot)) = \| \rho(t,\cdot) \|_{\dot{H}^{-1}(\T^d)} \geq \frac{ \| \rho(t,\cdot) \|^2_{L^2(\T^d)}}{ \| \rho(t,\cdot) \|_{\dot{H}^1(\T^d)} } 
\geq  \frac{ \| \bar\rho \|^2_{L^2(\T^d)}}{ \| \bar\rho \|_{\dot{H}^1(\T^d)} } \, e^{-Lt} \,.
\end{equation}
We observe that~\eqref{e:lowerfunctinter} has exactly the same time dependence as~\eqref{e:firstfct}. However, we derived~\eqref{e:lowerfunctinter} under stronger regularity assumptions on the initial datum; specifically, for~\eqref{e:lowerfunctinter} we require $\bar\rho \in \dot{H}^{1}(\T^d)$.

\begin{exercise}
Provide an alternative proof of~\eqref{e:propagPDE} by means of a direct energy estimate for~\eqref{e:PDE}, without relying on the representation formula~\eqref{e:character}. 
\end{exercise}

\begin{exercise}
Provide an alternative proof of~\eqref{e:lowerfunctinter} based on the duality formula
$$
\| f \|_{\dot{H}^{-1}(\T^d)} = \sup_{\| \varphi \|_{\dot{H}^{1}(\T^d)}=1} \int_{\T^d} f(x) \, \varphi(x) \, dx \,.
$$ 
In fact, this proof gives a better version of~\eqref{e:lowerfunctinter}, since it does not require any regularity of $\bar\rho$. 
\end{exercise}

\begin{remark}
Even though the DiPerna-Lions-Ambrosio~\cite{MR1022305,MR2096794} theory guarantees  well-posedness of~\eqref{e:PDE} below Lipschitz regularity, no analogue of the propagation of regularity~\eqref{e:propagPDE}, not even in a weaker form, holds in this regime, see the comments and references in Section~\ref{s:saturation}. The interpolation inequality~\eqref{e:interpol} is a true inequality: the two sides are not, in general, comparable.
\end{remark}

\subsection{Lower bound for the geometric mixing scale}\label{ss:easyflow}
We now present a direct proof of a lower bound for the geometric mixing scale, based on a purely Lagrangian argument. For simplicity, we restrict to the two-dimensional case and consider an initial configuration as in Figure~\ref{f:mixscale_1} (left), with value $+1$ on the grey region and $-1$ on the white region. We pick the point $\bar{x} = (\sfrac{3}{4} , \sfrac{1}{2})$ which is ``in the middle'' of the white region, 
and we denote ${\rm mix}_{\rm g} (\rho(t,\cdot)) = \varepsilon$ (recall the definition of the geometric mixing scale in Subsection~\ref{ss:defgeomix}). 


By the definition of geometric mixing scale, the ball of radius $\varepsilon$ centered at $\Phi(t,\bar{x})$ contains at least one point from the grey region. We denote $\Phi(t,\bar{y})$ such a point. Therefore,
$$
| \Phi(t,\bar{x}) - \Phi(t,\bar{y}) | \leq \varepsilon \,,
$$
and by the choice of the point $\bar{x}$
$$
|\bar{x}-\bar{y}| \geq \sfrac{1}{4} \,.
$$
As a consequence, we have
\begin{equation}\label{e:points}
\begin{aligned}
{\rm Lip}\, \Phi(t,\cdot)^{-1} &\geq
\frac{|\Phi(t,\cdot)^{-1}(\Phi(t,\bar{x})) - \Phi(t,\cdot)^{-1}(\Phi(t,\bar{y}))|}{|\Phi(t,\bar{x})-\Phi(t,\bar{y})|} \\
&= \frac{|\bar{x}-\bar{y}|}{|\Phi(t,\bar{x})-\Phi(t,\bar{y})|} \geq \frac{\sfrac{1}{4}}{\varepsilon} = \frac{1}{4\varepsilon} \,.
\end{aligned}
\end{equation}
By the assumption that ${\rm Lip}\, (u(t,\cdot)) \leq L$ for all $t\geq 0$, we have by~\eqref{e:regLipflowinverse} that
$$
{\rm Lip}\, \big( \Phi(t,\cdot)^{-1} \big) \leq e^{Lt} \,,
$$
which together with~\eqref{e:points} gives
$$
\varepsilon = {\rm mix}_{\rm g} (\rho(t,\cdot)) \geq \frac{1}{4} e^{-Lt} \,.
$$
This lower bound has exactly the same time dependence as~\eqref{e:firstfct} and~\eqref{e:lowerfunctinter}. The factor~$\sfrac{1}{4}$ encodes some (very weak) regularity of the initial configuration, namely the fact that there exists a ball with radius~$\sfrac{1}{4}$ entirely contained in the white region. 

\begin{remark}
The advantage of this geometric approach (compared with the energy estimates of Section~\ref{s:eneest} or the regularity–propagation estimates of Subsection~\ref{ss:regestinter}) is that it can be extended to settings of lower regularity, as we will carry out in Section~\ref{s:CDLproof}, based on the theory in~\cite{MR2369485}. 
\end{remark}

\section{Bressan's mixing scheme and a glimpse of DiPerna-Lions-Ambrosio theory}\label{s:bressan}

We have seen, by different proofs, that for a uniformly-in-time Lipschitz velocity field the (functional or geometric) mixing scale can decay at most at a negative exponential rate in time. In this section, we present a simple combinatorial example set on the two-dimensional torus (following the construction of~\cite{MR2033002}, see also~\cite{MR2009116,MR3026556}) which shows, via a ``slice-and-dice'' mechanism, that both mixing scales can indeed decrease to zero at a negative exponential rate. However, the velocity field in this example is not Lipschitz: it is built from piecewise-constant shear flows and has discontinuities. Still, in a lower regularity setting, the velocity enjoys a uniform-in-time bound on its first-order spatial derivative. In fact, it is uniformly bounded in time in $BV$ (bounded variation). We will not develop the full rigorous theory of $BV$ spaces here, but only provide some intuition: in our case the velocity field is piecewise constant, with jumps along lines. In a sense, the example in this section shows that exponential mixing is sharp once we relax the regularity class assumed for the velocity field.

To describe the evolution of the passive scalar, we first specify the configurations obtained at discrete times. We then show that there exists a velocity field for which these prescribed configurations coincide with the solutions of~\eqref{e:PDE} at the given discrete times. 

\begin{figure}[h]
\begin{center}
\includegraphics[scale=0.25]{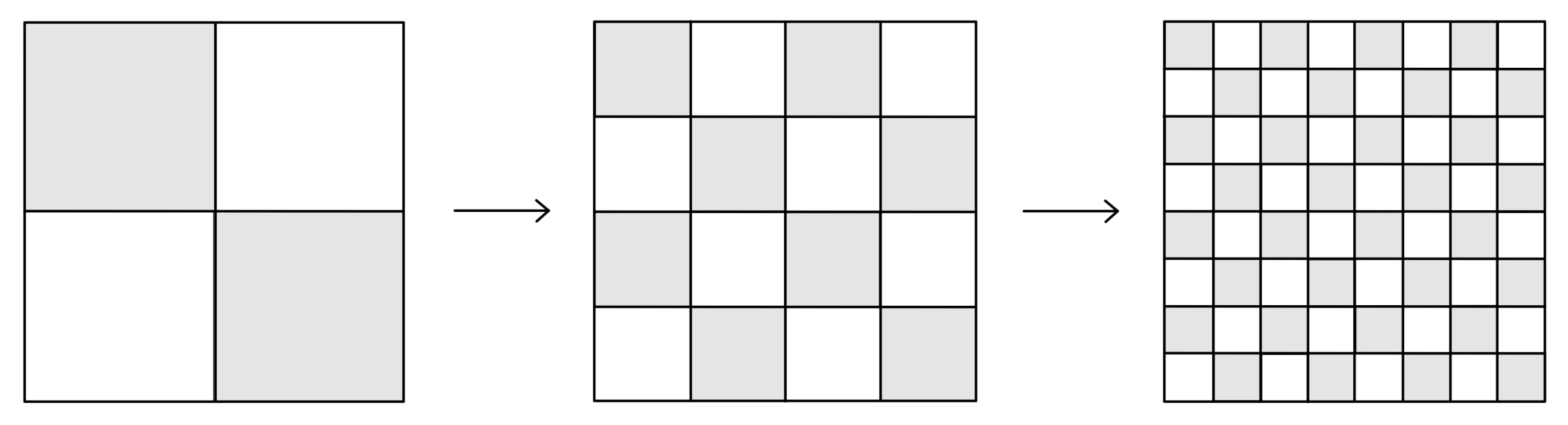}
\end{center}
\caption{Bressan's mixing scheme}\label{f:bresteps}
\end{figure}

The configuration at the dyadic times $t_0$, $t_1$, $t_2$, \ldots, with $t_k = \sum_{j=1}^k 2^{-j}$, is depicted in Figure~\ref{f:bresteps}. At the initial time $t_0=0$, the configuration takes the alternating values $\pm 1$ on a checkerboard with squares of side $\sfrac{1}{2}$. At each step, the side of the squares is reduced by a factor of $\sfrac{1}{2}$, so that at the discrete times~$t_k$ the side decreases dyadically as $2^{-(k+1)}$, with the solution taking the values $\pm 1$ on the corresponding squares. At the final time $t = 1$, the solution is perfectly mixed: the configurations converge weakly to zero as the side length of the squares tends to zero.

\begin{figure}[h]
\begin{center}
\includegraphics[scale=0.25]{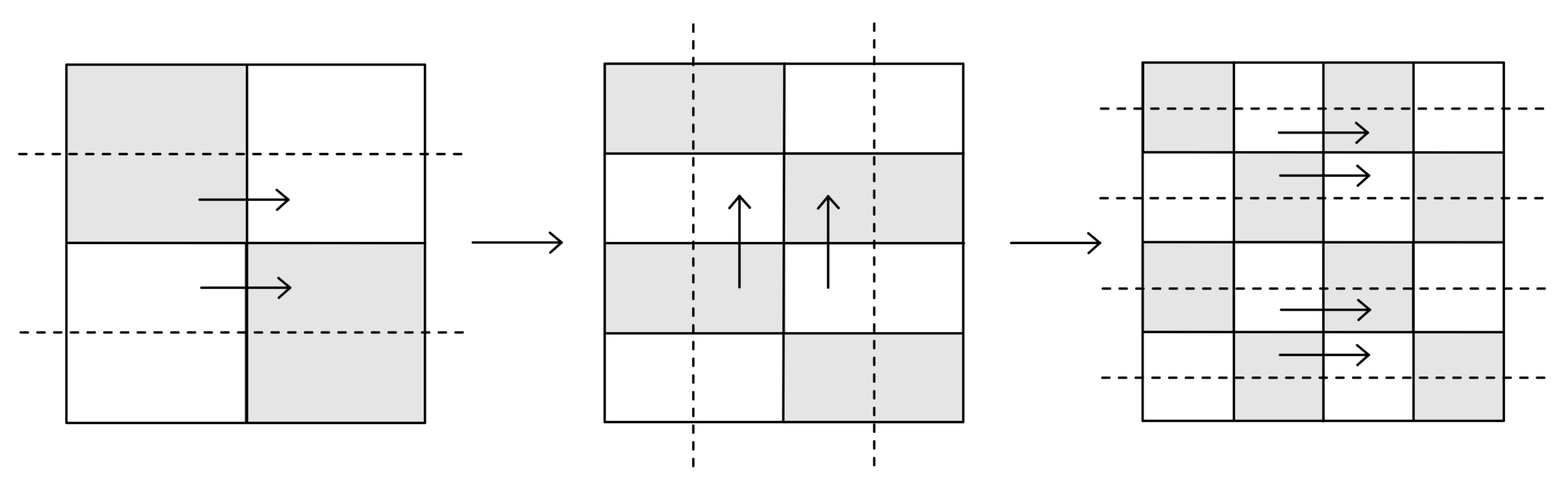}
\caption{Intermediate stages and shear flows in Bressan's mixing scheme.}\label{f:brestages} 
\end{center}
\end{figure}

This evolution of the passive scalar configuration can be realized by a velocity field composed of piecewise constant shear flows, as can be seen from Figure~\ref{f:brestages}. Within each step of the dyadic evolution, we insert an intermediate stage: first, the squares are split in one direction and become rectangles; second, the rectangles are split again, producing squares of the prescribed size, which are  arranged in the correct positions of the finer checkerboard, by the periodicity of the torus. As the checkerboard becomes finer and finer, the corresponding shear flows also become increasingly refined. Since the velocity field is parallel to the discontinuity interfaces, it is easy to check that the velocity field is divergence-free.

Let us now carry out a scaling analysis to understand this evolution quantitatively. At the $k$-th step, the number of discontinuity interfaces is of the same order as the number of channels, namely $2^k$, so the total length of the interfaces is also of order $2^k$. The displacement is of order $2^{-k}$. With the chosen time step~$2^{-k}$, this requires a velocity of order $1$ inside the active channels, and $0$ elsewhere. Thus, the velocity field is uniformly-in-time bounded in $L^\infty$. For a piecewise constant velocity field of this type, the $BV$-norm (that is, rigorously, the norm of the distributional derivative of first order in the sense of measures) is given by the product of the size of the jumps of the velocity field and the total length of the interfaces. Hence, at the $k$-th step, it is of order $2^k$. Rewriting this as a function of time, we obtain that $\| u(t,\cdot) \|_{\dot{BV}(\T^2)} \sim \frac{1}{1-t}$. 

Notice that, by time reversal, we have produced an example of nonuniqueness for the continuity equation~\eqref{e:PDE}. This was expected, based on the linear-in-time lower bound established in Subsection~\ref{ss:kinPDE} under uniform-in-time integrability bounds on the velocity field. In fact, observing that $t \mapsto \frac{1}{1-t}$ is not integrable at $t = 1$, this construction provides an instance of the sharpness of the DiPerna-Lions-Ambrosio~\cite{MR1022305,MR2096794} theory (see also the review paper~\cite{MR3283066}). Indeed, this theory asserts (to keep things simple, for a bounded, divergence-free velocity field on the torus, even though these assumptions can be relaxed) that, assuming
$$
u \in L^1([0,T] ; W^{1,p}(\T^d)) \;\; \text{for some $1 \leq p \leq \infty$}
\quad \text{ or } \quad
u \in L^1([0,T] ; BV(\T^d)) \,,
$$
then the continuity equation~\eqref{e:PDE} admits a unique bounded weak solution, and there exists a unique suitable selection of solutions of~\eqref{e:ODE} giving rise to the so-called regular Lagrangian flow associated with the velocity field. Moreover, stability holds in both cases. We will describe more refined examples showing the sharpness of these assumptions in Section~\ref{s:saturation}. The DiPerna-Lions-Ambrosio~\cite{MR1022305,MR2096794} theory is based on the notion of renormalized solution and on the study of the convergence of suitable commutators, and it has been extended in~\cite{MR4071413} to the case of nearly incompressible $BV$ velocity fields, also answering in the positive the so-called Bressan's compactness conjecture (this conjecture asserts the compactness of sequences of regular Lagrangian flows under incompressibility and $BV$ bounds, see~\cite{MR2033003}). 

We now turn to the time asymptotics of the decay of the mixing scale in this construction, introducing a rescaling of time. Instead of dyadic time steps, we take unit time steps, that is, we choose $\tilde{t}_k = k$. At any $\tilde{t}_k$, the configuration coincides with that at $t_k$ in the original evolution: we are just ``delaying time''. In this way, the weak convergence of the passive scalar to zero occurs asymptotically as $t \to \infty$. In the corresponding scaling analysis, the velocity field now acts over a time of order $1$ to produce a displacement of order $2^{-k}$. Hence the velocity takes values $2^{-k}$ or $0$. The total length of discontinuity interfaces does not change, so the $BV$-norm is uniformly bounded in time. Meanwhile, the mixing scale of the passive scalar decays exponentially in time. Thus, after this rescaling we obtain
$$
\| u(t,\cdot) \|_{\dot{BV}(\T^2)} \sim 1
\qquad \text{ and } \qquad
{\rm mix}_{\rm g} (\rho(t,\cdot)) \sim {\rm mix}_{\rm f} (\rho(t,\cdot)) \sim 2^{-t} \,,
$$
providing an example of exponential-in-time decay of the mixing scale.

\section{An integral functional measuring the logarithmic regularity of the flow}\label{s:functplan}

The linear-in-time lower bound on the mixing scale obtained in Subsection~\ref{ss:ensPDE} leaves open the possibility of perfect mixing in finite time when the velocity field has uniformly-in-time bounded enstrophy. However, this would contradict the results of the DiPerna-Lions-Ambrosio~\cite{MR1022305,MR2096794} theory, as recalled at the end of the previous section. 

Our next goal  is to establish sharper lower bounds on the mixing scale in the case of velocity fields with uniformly-in-time bounded enstrophy. Rather than proceeding by energy estimates, we follow the approach in~\cite{MR2369485} (see also the review paper~\cite{MR3283066}) and extend the Lagrangian argument from Subsection~\ref{ss:easyflow}. The theory developed in~\cite{MR2369485} is based on estimates for the regular Lagrangian flow and is, by design, fully quantitative. This stands in contrast to the DiPerna-Lions-Ambrosio~\cite{MR1022305,MR2096794} theory, which is non-quantitative, relying instead on the concept of renormalized solutions and abstract compactness arguments. In place of a Gr\"onwall estimate, we will rely on another formal estimate for the regularity of the flow, namely
\begin{equation}\label{e:formallog}
\begin{aligned}
\frac{d}{dt} \log | \nabla \Phi| &\leq
\frac{|\nabla \dot \Phi|}{|\nabla \Phi|} =
\frac{|\nabla (u(\Phi))|}{|\nabla\Phi|} \leq
\frac{|\nabla u|(\Phi) \, |\nabla\Phi|}{|\nabla\Phi|} = |\nabla u|(\Phi) \,,
\end{aligned}
\end{equation}
which implies
$$
\log |\nabla \Phi| \lesssim \int_0^t |\nabla u|(\Phi) \, ds \,.
$$
This formally means that the logarithm of the derivative of the flow can be controlled by the derivative of the velocity field along the flow. If the velocity field is assumed to be Lipschitz in the spatial variable, this estimate is equivalent to the classical Gr\"onwall inequality. 

However, the logarithmic estimate~\eqref{e:formallog} has the advantage of allowing to obtain results for velocity fields with Sobolev regularity in the spatial variable, as first shown in~\cite{MR2369485}. In order to do so, we define an integral functional depending on the flow:
\begin{equation}\label{e:intfunc}
\mathcal{G}(\Phi) =
\left\| \sup_{0 \leq t \leq T} \sup_{r>0} 
\aint_{B(x,r)} \log \left(1+ \frac{|\Phi(t,x)-\Phi(t,y)|}{r} \right) \, dy 
\right\|_{L^p_x(\T^d)} \,.
\end{equation}
The functional defined in~\eqref{e:intfunc} is a discrete, integral version of the quantity estimated in~\eqref{e:formallog}. There are three main steps in the argument to show exponential-in-time lower bounds on the mixing scale for Sobolev velocity fields, under the assumption
$$
u \in L^\infty([0,+\infty[; W^{1,p}(\T^d)) \;\; \text{for some $1 < p \leq \infty$}
$$
(notice that the endpoint case $p=1$ is not covered, and therefore not $BV$ either):
\begin{itemize}
\item Establish bounds on $\mathcal{G}(\Phi)$ in term of norms of the velocity field. Here (and only here) we use the dynamics, that is, the fact that~$\Phi$ is a solution of~\eqref{e:ODE} (see Subsection~\ref{s:quaest}). 
\item Use the bounds on $\mathcal{G}(\Phi)$ to show some mild regularity of the flow, more precisely, a quantitative Lusin-Lipschitz estimate (see Subsection~\ref{s:lusin}). 
\item Show that the Lusin-Lipschitz estimate implies exponential-in-time lower bounds on the mixing scale (see Subsection~\ref{s:lower}).
\end{itemize}

First of all, we differentiate in time the integral in the definition of the functional $\mathcal{G}(\Phi)$ in~\eqref{e:intfunc}, obtaining
\begin{equation}\label{e:firstcomput}
\begin{aligned}
\frac{d}{dt} \aint_{B(x,r)} \log \left(1+ \frac{|\Phi(t,x)-\Phi(t,y)|}{r} \right) \, dy 
&\leq 
\aint_{B(x,r)} \frac{|u(t,\Phi(t,x))-u(t,\Phi(t,y))|}{r+|\Phi(t,x)-\Phi(t,y)|} \, dy \\
&\leq
\aint_{B(x,r)} \frac{|u(t,\Phi(t,x))-u(t,\Phi(t,y))|}{|\Phi(t,x)-\Phi(t,y)|} \, dy \,.
\end{aligned}
\end{equation}
If the velocity field has spatial Lipschitz regularity, the right hand side of the last expression can be easily bounded by a constant. In order to obtain bounds on the difference quotients in the case of Sobolev velocity fields we need some tools from harmonic analysis, that we introduce in the next section.

\section{Some tools from harmonic analysis}\label{s:hatools}

We know thanks to the Lebesgue differentiation theorem that  for $f \in L^1(\T^d)$ it holds
$$
\aint_{B(x,r)} f(y) \, dy \to f(x) 
\qquad\text{as $r\to 0$, for a.e.~$x\in\T^d$.}
$$
Instead of the limit as $r \to 0$, let us consider the associated maximal operator:
$$
Mf(x) = \sup_{r>0} \aint_{B(x,r)} |f(y)| \, dy 
\qquad \text{ for~$x\in\T^d$.}
$$
\begin{exercise}
Write the definition of the maximal function of a measure.
\end{exercise}
If $f\in L^1(\T^d)$, then $Mf(x)$ is finite for a.e.~$x\in\T^d$ and it is called the maximal function of $f$. This defines an operator which can be shown to enjoy the following strong estimates:
\begin{equation}\label{e:strong}
\| Mf \|_{L^p(\T^d)} \leq C_{d,p} \| f \|_{L^p(\T^d)}
\qquad \text{for $1<p\leq \infty$.}
\end{equation}

\begin{exercise} 
Show the strong estimate for $p=\infty$. Understand why the strong estimate cannot hold for $p=1$. 
\end{exercise}

For $p=1$ only a weak estimate holds, namely
\begin{equation}\label{e:weak}
||| Mf |||_{L^{1,\infty}(\T^d)} \leq C_{d,1} \|f\|_{L^1(\T^d)} \,,
\end{equation}
where for any measurable function $f$ defined on $\T^d$ we set
\begin{equation}\label{e:weaknorm}
||| f |||_{L^{1,\infty}(\T^d)} = \sup_{\lambda>0} \Big\{ \lambda \big| \{ x \in \T^d \;:\; |f(x)|>\lambda \} \big| \Big\} \,.
\end{equation}
The space $L^{1,\infty}(\T^d)$ is defined as the space of all measurable functions such that the quantity defined in~\eqref{e:weaknorm} is finite. 

\begin{exercise}
Use Chebyshev’s inequality to show that $L^1(\T^d) \subset L^{1,\infty}(\T^d)$. Show that the inclusion is proper, that is, there exists $f \in L^{1,\infty}(\T^d) \setminus L^1(\T^d)$. 
\end{exercise}

\begin{remark}
The quantity defined in~\eqref{e:weaknorm} does not define a norm (it is not subadditive). For this reason, we have chosen the notation with the triple bar.
\end{remark}

\begin{remark}
The weak estimate~\eqref{e:weak} can be proven using Vitali's covering lemma (see Section~\ref{s:lower} below for its statement). The strong estimate~\eqref{e:strong} for $1<p<\infty$ follows from the weak estimate~\eqref{e:weak} and the strong estimate~\eqref{e:strong} for $p=\infty$ by an interpolation argument.
\end{remark}

The necessity of introducing maximal functions in our context arises from their usefulness in estimating difference quotients. In fact, it can be shown that for every $f \in W^{1,1}(\T^d)$ (or even $f \in BV(\T^d)$) there exists a negligible set $N \subset \T^d$ such that
\begin{equation}\label{e:increments}
| f(x) - f(y) | \leq C_d \, |x-y| \, \big( M Df (x) + M Df (y) \big)
\qquad \forall \, x,y \in \T^d \setminus N \,.
\end{equation} 
In the case of a function $f$ with spatial Lipschitz regularity, it holds $Df \in L^\infty(\T^d)$ and therefore also $M Df \in L^\infty(\T^d)$. Therefore,~\eqref{e:increments} expresses the usual estimate for the difference quotients of a Lipschitz function. In the case of a Sobolev or $BV$ function,~\eqref{e:increments} expresses a bound on the difference quotients which is not uniform, but rather depends on on the size of the maximal function. Very informally,~\eqref{e:increments} can be interpreted as a Lipschitz estimate with a point-dependent constant.


\begin{exercise}
Prove the estimate~\eqref{e:increments} exploiting the following idea.
 Fix $x,y \in \T^d$. If~$f$ is Lipschitz, we can estimate $|f(x)-f(y)|$ by the integral of $Df$ over the segment connecting $x$ and $y$. But if $f$ is only Sobolev, such an integral might be infinite. Set $r=|x-y|$ and define $C_{x,y} = B(x,r) \cap B(y,r)$. For every $z \in C_{x,y}$, estimate $|f(x)-f(y)|$ by the integral of $Df$ over the segment connecting $x$ and $z$, plus the integral of $Df$ over the segment connecting $z$ and $y$. Average over $z \in C_{x,y}$ and employ some changes of variable to obtain the maximal function of $Df$. In the proof, you can assume $f$ to be smooth, but it is important that the estimate you obtain only depends on the Sobolev or $BV$ norm of $f$. 
\end{exercise}

\section{Estimates for uniformly-in-time Sobolev velocity fields}\label{s:CDLproof}

In this section, we carry out the plan described in Section~\ref{s:functplan} in order to establish an exponential-in-time lower bound on the mixing scale in the case of uniformly-in-time Sobolev velocity fields, following the approach developed in~\cite{MR2369485}. The proof is divided into the following three subsections.

\subsection{Quantitative estimates for the functional $\mathcal{G}(\Phi)$}\label{s:quaest}
Proceeding from the computation in~\eqref{e:firstcomput} and using~\eqref{e:increments}, we obtain
$$
\begin{aligned}
\frac{d}{dt} \aint_{B(x,r)} & \log \left(1+ \frac{|\Phi(t,x)-\Phi(t,y)|}{r} \right) \, dy 
\leq 
\aint_{B(x,r)} \frac{|u(t,\Phi(t,x))-u(t,\Phi(t,y))|}{|\Phi(t,x)-\Phi(t,y)|} \, dy \\
&\lesssim
\aint_{B(x,r)} \left[ M Du (t,\cdot) (\Phi(t,x))  + M Du (t,\cdot) (\Phi(t,y)) \right] \, dy \\
&= M Du (t,\cdot) (\Phi(t,x))
+
\aint_{B(x,r)}  M Du (t,\cdot) (\Phi(t,y))  \, dy \,.
\end{aligned}
$$
Therefore, it follows
$$
\begin{aligned}
\aint_{B(x,r)} & \log \left(1+ \frac{|\Phi(t,x)-\Phi(t,y)|}{r} \right) \, dy 
\lesssim \aint_{B(x,r)} \log \left(1+ \frac{|x-y|}{r}\right) \, dy \\
& + \int_0^t M Du (s,\cdot) (\Phi(s,x)) \, ds
+
\int_0^t \aint_{B(x,r)}  M Du (s,\cdot) (\Phi(s,y))  \, dy \, ds \,.
\end{aligned}
$$
Taking the supremum over $0 \leq t \leq T$ and $r>0$ and the $L^p$-norm in the variable $x$, we get
$$
\begin{aligned}
\mathcal{G}(\Phi) 
&\lesssim
1 +
\int_0^T \| M Du (s,\cdot) (\Phi(s,x))\|_{L^p_x(\T^d)} \, ds
+
\int_0^T \| M \left[    M Du (s,\cdot) (\Phi(s,\cdot)) \right] (x) \|_{L^p_x(\T^d)} \, ds \\
&\lesssim 1+\int_0^T \| M Du (s,\cdot) (\Phi(s,x))\|_{L^p_x(\T^d)} \, ds \\
&= 1+\int_0^T \| M Du (s,\cdot) (x)\|_{L^p_x(\T^d)} \, ds \\
&\lesssim 1+\int_0^T \| Du (s,\cdot) (x)\|_{L^p_x(\T^d)} \, ds = 1+ \| Du\|_{L^1([0,T];L^p(\T^d))} \,.
\end{aligned}
$$
We conclude that there exists $\bar{C} = \bar{C} (d,p,T)$ such that
\begin{equation}\label{e:Gquant}
\mathcal{G}(\Phi) \leq \bar{C} \big( 1+\| Du\|_{L^1([0,T];L^p(\T^d))} \big) \,.
\end{equation}
In an informal sense, estimate~\eqref{e:Gquant} indicates that the regular Lagrangian flow possesses a logarithmic derivative. Alternative approaches involving logarithmic function spaces have been developed in~\cite{MR4263701,MR4377866}, removing the divergence-free condition in~\cite{MR3862947}, and through Littlewood–Paley theory in~\cite{MR4747245}.

\subsection{Quantitative Lusin-Lipschitz regularity of the flow}\label{s:lusin}
We now show how a mild regularity property of the regular Lagrangian flow follows from the bound on $\mathcal{G}(\Phi)$. Fix $\eta>0$.  Using Chebyshev’s inequality applied to
$$
\mathcal{G}(\Phi)  = \left\| \; \ldots \; \right\|_{L^p_x}
\qquad \text{ and } \qquad
\lambda = \frac{\bar{C} (1+\| Du\|_{L^1_t L^p_x})}{\eta^{\sfrac{1}{p}}} \,,
$$
we deduce the existence of a compact set $K \subset \T^d$ with $|\T^d \setminus K| \leq \eta$ such that
\begin{equation}\label{e:aftercheb}
\sup_{0 \leq t \leq T} \sup_{r>0} 
\aint_{B(x,r)} \log \left(1+ \frac{|\Phi(t,x)-\Phi(t,y)|}{r} \right) \, dy 
\leq 
\frac{\bar{C} (1+\| Du\|_{L^1_t L^p_x})}{\eta^{\sfrac{1}{p}}}
\qquad
\forall \, x \in K \,.
\end{equation}
We now show that $\Phi(t,\cdot)$ is Lipschitz on the set $K$, with a quantitative estimate on its Lipschitz constant. 
From~\eqref{e:aftercheb}, we deduce
$$
\aint_{B(x,r)} \log \left(1+ \frac{|\Phi(t,x)-\Phi(t,y)|}{r} \right) \, dy 
\leq 
\frac{\bar{C} (1+ \| Du\|_{L^1_t L^p_x})}{\eta^{\sfrac{1}{p}}}
\qquad
\forall \, 0 \leq t \leq T\,, \; r>0\,,  \; x \in K \,.
$$
Fix $x,x' \in K$, set $r=|x-x'|$, and define  $C_{x,x'} = B(x,r) \cap B(x',r)$. 

\begin{exercise}
Show that $| B(x,r) | = c_d |C_{x,x'}| $ for some dimensional constant $c_d$. 
\end{exercise}

\begin{exercise}
The function $z\mapsto\log(1+z)$ defined for $z\geq 0$ is nonnegative, increasing, and subadditive, that is
$$
\log (1+x+y) \leq \log(1+x) + \log(1+y)
\qquad \forall \, x,y \geq 0 \,.
$$
\end{exercise}

Using the last two exercises we deduce that
$$
\begin{aligned}
\log&\left( 1 + \frac{|\Phi(t,x)-\Phi(t,x')|}{r}\right) 
=
\aint_{C_{x,x'}} \log\left( 1 + \frac{|\Phi(t,x)-\Phi(t,x')|}{r}\right)  \, dy \\
&\leq 
\aint_{C_{x,x'}} \log\left( 1 + \frac{|\Phi(t,x)-\Phi(t,y)|}{r}\right)  \, dy 
+
\aint_{C_{x,x'}} \log\left( 1 + \frac{|\Phi(t,x')-\Phi(t,y)|}{r}\right)  \, dy \\
&\leq 
c_d \; \aint_{B(x,r)} \log\left( 1 + \frac{|\Phi(t,x)-\Phi(t,y)|}{r}\right)  \, dy 
+
c_d \; \aint_{B(x',r)} \log\left( 1 + \frac{|\Phi(t,x')-\Phi(t,y)|}{r}\right)  \, dy \\
&\leq 
2 c_d \frac{\bar{C} (1+\| Du\|_{L^1_t L^p_x})}{\eta^{\sfrac{1}{p}}} \,.
\end{aligned}
$$
We have thus shown the following: for any $\eta>0$ there is $K \subset \T^d$ (which does not depend on time) with~$| \T^d \setminus K| \leq \eta$ such that
\begin{equation}\label{e:luslip}
| \Phi(t,x) - \Phi(t,x')| \leq |x-x'| \exp \left( 2 c_d \frac{\bar{C} (1+\| Du\|_{L^1_t L^p_x})}{\eta^{\sfrac{1}{p}}} \right)
\qquad \text{ for all $x,x'\in K$ and $0 \leq t \leq T$. }
\end{equation}
This is what we call quantitative Lusin-Lipschitz regularity of the regular Lagrangian flow.

\begin{remark}
In general, if the velocity field $u$ is not Lipschitz, the regular Lagrangian flow $\Phi$ is not Lipschitz. In fact, it might even be discontinuous. 
\end{remark}

\begin{exercise}
Understand how the Lusin-Lipschitz estimate~\eqref{e:luslip} implies strong compactness of regular Lagrangian flows under suitable bounds, and understand how this can be exploited to prove existence of the regular Lagrangian flow. 
\end{exercise}

\begin{remark}
The approach developed in the last two subsections can also be applied to establish uniqueness and stability of the regular Lagrangian flow under suitable assumptions and can be seen as a quantitative version of the DiPerna-Lions-Ambrosio~\cite{MR1022305,MR2096794} theory. It can be carried out for velocity fields in $W^{1,1}$ by means of an interpolation argument exploiting the factor $r$ in the denominator of~\eqref{e:firstcomput}, see~\cite{MR3207161}. This extension, however, concerns only the uniqueness and stability results, not the Lusin-Lipschitz estimate. The case of $BV$ velocity fields has been treated in~\cite{MR4242824}, where the resulting estimates involve the polar decomposition of the singular part of the derivative of the velocity field.
\end{remark}

\begin{remark}\label{r:backflow}
The Lusin-Lipschitz estimate~\eqref{e:luslip} can also be established for the two-parameter regular Lagrangian flow $\Phi = \Phi(s,t,x)$ solving
$$
\left\{ \begin{array}{l}
\partial_t \Phi (s,t,x) = u (t,\Phi(s,t,x)) \\
\Phi(s,x) = x  
\end{array}\right.
$$
and, in particular, it holds for the inverse of the regular Lagrangian flow as well.
\end{remark}

\begin{remark}
The above estimates in the case of the space $\R^d$ can be localized. Moreover, they can also be extended to some cases when the divergence-free condition is dropped. 
\end{remark}

\subsection{Lower bounds on the mixing scale for uniformly-in-time Sobolev velocity fields}\label{s:lower}

The idea is to follow the same approach as in Subsection~\ref{ss:easyflow}, but now exploiting the fact that the flow $\Phi$ is Lipschitz on large subsets, together with a suitable covering argument.

We focus on the two-dimensional case of $\T^2$ and we consider the initial configuration on the left in Figure~\ref{f:mixscale_1}. 
We denote the left half of the torus by $A$, and the right half by $A'$. Also, we fix $t>0$, let $\Phi = \Phi(t,\cdot)$, and denote by $\varepsilon$ the geometric mixing scale at time $t$.

Thanks to the Lusin-Lipschitz estimate~\eqref{e:luslip}, together with Remark~\ref{r:backflow}, given $\eta>0$ there is a set~$H \subset \T^2$ with $|H| \leq \eta$ such that
$$
{\rm Lip}\, \left( \Phi^{-1} |_{\T^2\setminus H} \right) \leq \exp(\beta t M) \,,
$$
for some constant 
$$
\beta = \beta(2,p,\eta,T)
$$
and where 
$$
M = 1+ \sup_{0 \leq s \leq t} \| \nabla u(s,\cdot)\|_{L^p(\T^2)} \,.
$$
By the assumption on the geometric mixing scale at time $t$ we have
\begin{equation}\label{e:A}
| B_\eps(\Phi(x)) \cap \Phi(A')|
\geq
\kappa | B_\eps(\Phi(x))|
\qquad \forall \, x \in A\,.
\end{equation}
Let us define
\begin{equation}\label{e:tildeA}
\tilde{A} = \big\{ x\in A \;:\; B_\eps(\Phi(x)) \cap [ \Phi(A') \setminus H ] = \emptyset \big\} \,.
\end{equation}
\begin{exercise}
Check the implication
\begin{equation}\label{e:sets}
V_1 \cap (V_2\setminus V_3) = \emptyset
\qquad \Longrightarrow \qquad
V_1\cap V_2 \subset V_1 \cap V_3
\end{equation}
for any triple of sets $V_1$, $V_2$, and $V_3$. 
\end{exercise}
Implication~\eqref{e:sets} allows to deduce from~\eqref{e:A} and~\eqref{e:tildeA} that
\begin{equation}\label{e:AtildeA}
| B_\eps(\Phi(x)) \cap H| \geq \kappa | B_\eps(\Phi(x))|
\qquad \forall \, x \in \tilde{A} \,.
\end{equation}

We now recall Vitali's covering lemma: {\sl Let $\mathcal{F}$ be a family of balls with uniformly bounded radius contained in the torus. Then, there exists a countable, disjoint subfamily $\mathcal{F}'$ such that
$$
\bigcup_{B \in \mathcal{F}} B
\subset
\bigcup_{B \in \mathcal{F}'} 5B \,,
$$
where, for a ball $B$, we denote by $5B$ the ball with the same center as $B$ and radius five times larger.}

Exploiting Vitali's covering lemma and~\eqref{e:AtildeA} we deduce that
$$
| \Phi(\tilde{A})| \leq \frac{5^2}{\kappa} |H| \leq \frac{5^2}{\kappa} \eta \,.
$$
By the incompressibility if the flow this implies that
$$
| \tilde{A} | \leq \frac{5^2}{\kappa} \eta \,.
$$
Again incompressibility gives
$$
| \Phi^{-1} (H)| = |H| \leq \eta \,.
$$
We now choose $\eta = \eta(2,\kappa)$ (consequently, $\beta=\beta(2,p,\kappa,T)$) so that
$$
| \tilde{A}| + | \Phi^{-1}(H )| < \sfrac{1}{6} \,.
$$
Then, there exists 
$$
\bar{x} \in A \setminus \left( \tilde{A} \cup \Phi^{-1}(H) \right)
\qquad \text{ such that } \qquad
{\rm dist}\, (\bar{x} , A') \geq \sfrac{1}{6} \,.
$$


Denote $\bar{y} = \Phi(\bar{x})$. Since $\bar{x} \not \in \tilde{A}$, there exists $\bar{z} \in B_\eps(\bar{y}) \cap [ \Phi(A') \setminus H]$. Clearly
$$
|\bar{y}-\bar{z}| \leq \eps \,,
\qquad
|\bar{x}-\Phi^{-1}(\bar{z})| \geq \sfrac{1}{6} \,,
\qquad \text{ and } \qquad
\Phi^{-1}(\bar{z}) \in A' \,.
$$
On the other hand, $\bar{y},\bar{z} \not \in H$. 
Therefore
$$
\frac{\sfrac{1}{6}}{\eps} \leq {\rm Lip} \, \left( \Phi^{-1} |_{\T^2\setminus B} \right) \leq \exp( \beta t M) \,,
$$
which implies
$$
\eps \geq \frac{1}{6} \exp( -\beta t M) \,,
$$
which recalling our notations means
$$
{\rm mix}_{\rm g} \, (\rho(t,\cdot) \geq \frac{1}{6} \exp \left( - \beta \big(1+\| Du \|_{L^\infty_t L^p_x}\big) \, t \right) \,.
$$
Thus, we have established a lower bound that decays exponentially in time (the same rate as in the Lipschitz case, recall Subsection~\ref{ss:lipPDE} and Section~\ref{s:flowestimates}) for velocity fields that are Sobolev in space, more precisely in $L^\infty([0,T] ; W^{1,p}(\T^d))$ for some $1 < p \leq \infty$. In particular, this improves substantially upon the linear-in-time lower bound proved in Subsection~\ref{ss:ensPDE}, and it rules out perfect mixing in finite time, in accordance with the DiPerna–Lions theory.

This result, however, is limited to $p > 1$, due to the lack of a strong inequality for the maximal function (recall Section~\ref{s:hatools}), or alternatively, due to the lack of convexity of the quantity $f \mapsto ||| f |||_{L^{1,\infty}(\T^d)}$. It remains an open question whether the same exponential lower bound holds for $BV$ velocity fields, a problem known as the Bressan mixing conjecture (see~\cite{bressan2006prize}). Lower bounds on the mixing scale with concatenated exponential behavior have been established in~\cite{2504.03023}.

Other approaches based more on harmonic analysis can be found in~\cite{MR3799259,2402.11642}. The passage from the geometric to the functional mixing scale is highly nontrivial: this  has been  established in~\cite{MR3141856,MR3207161}. Stronger lower bounds can be obtained in cases with specific geometric structure, such as Hamiltonian flows~\cite{MR4797689} or cellular flows~\cite{MR3703559}.

\section{Optimality of the estimates and properties of the regular Lagrangian flow}\label{s:saturation}

Let us summarize the results obtained so far. Considering either the geometric or functional mixing scale, which we now denoted generically by ${\rm mix}(\rho)$ in what follows, we have shown that:
\begin{itemize}
\item If $u \in L^\infty([0,+\infty[ ; W^{1,p}(\T^d))$ for some $1 < p \leq \infty$, then for any initial configuration with some regularity it holds ${\rm mix}(\rho(t,\cdot)) \gtrsim e^{-ct}$ for all $t\geq 0$.
\item There exist a velocity field $u \in L^\infty([0,+\infty[ ; BV(\T^d))$ and an initial configuration such that it holds 
${\rm mix}(\rho(t,\cdot)) \simeq e^{-ct}$ for all $t\geq 0$.
\end{itemize}
At the same time, Bressan's mixing conjecture~\cite{bressan2006prize} is still open:
\begin{itemize}
\item Is it true that if $u \in L^\infty([0,+\infty[; BV(\T^d))$, then for any sufficiently regular initial configuration, one has ${\rm mix}(\rho(t,\cdot)) \gtrsim e^{-ct}$ for all $t\geq 0$?
\end{itemize}
As already noted, extending the approach in~\cite{MR2369485} to the case $p = 1$ appears to be very difficult, since this is an endpoint for the harmonic analysis techniques employed there.
In a somewhat opposite direction, in~\cite{MR3904158} (see also the announcement~\cite{MR3268760}) the optimality of the estimates and rates  in~\cite{MR2369485} has been established, by constructing examples that can be viewed, informally, as extensions of the combinatorial construction in Section~\ref{s:bressan} to the setting of Sobolev-regular, or even Lipschitz, velocity fields.
In this section, we first briefly recall these constructions and then discuss what they reveal about the geometry and regularity of regular Lagrangian flows associated with Sobolev velocity fields.

\begin{figure}[h]
\begin{center}
\includegraphics[scale=0.5]{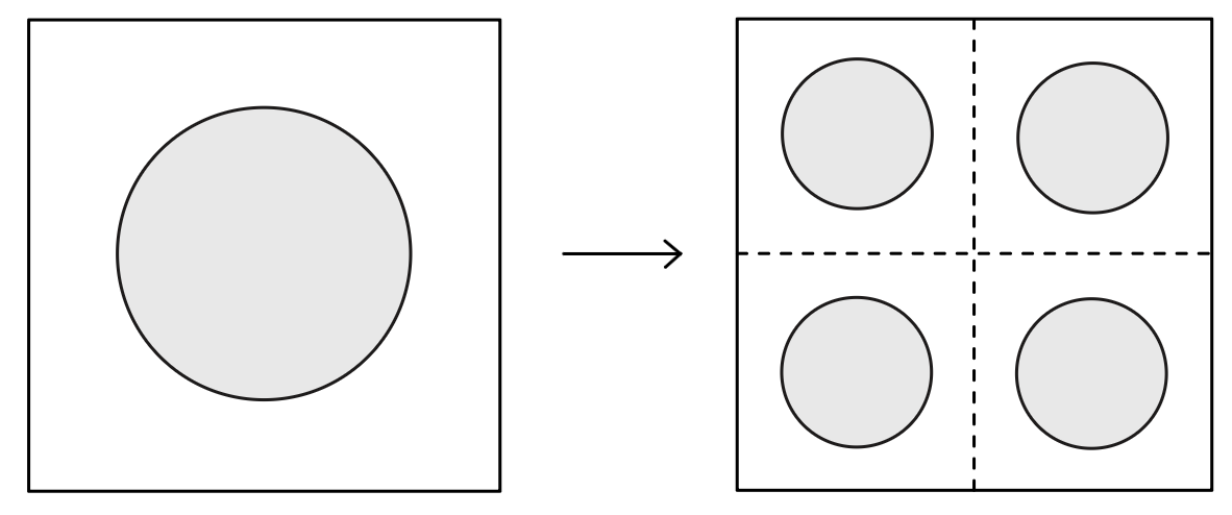}
\caption{The first step of a self-similar evolution.}\label{f:disks}
\end{center}
\end{figure}

The constructions in~\cite{MR3904158} share with the example in Section~\ref{s:bressan} the property of being self-similar: at each integer time, the configuration is a rescaled copy of a fixed portion of the configuration at the previous integer time.
However, unlike the example in Section~\ref{s:bressan}, these configurations are no longer based on a checkerboard pattern, and simply regularizing the shear flows from that example is not sufficient to obtain a velocity field with Sobolev regularity instead of merely $BV$. Instead, one must consider a geometric shape such as a disk, and construct a genuinely two-dimensional velocity field which at some finite time transforms the disk into four disjoint disks of equal size and the same total area as the initial disk, as in Figure~\ref{f:disks}. 
Inspecting the evolution shown in Figure~\ref{f:disks}, we immediately observe that it cannot be generated by a Lipschitz velocity field. Indeed, the flow associated with a Lipschitz velocity field is itself Lipschitz in the spatial variable and, in particular, a homeomorphism. Consequently, it cannot alter the topology of sets. In Figure~\ref{f:disks}, however, we see that a connected set evolves into a disconnected one.
As shown in~\cite{MR3904158}, using a technique based on deformation and pinching of sets, such an evolution can nevertheless be realized as the flow of a velocity field belonging to any fractional Sobolev space in the spatial variable $L^\infty([0,T];W^{s,p}(\T^2))$, provided that $s$ and $p$ lie below the threshold for embedding into the class of Lipschitz velocity fields,  in particular, for any Sobolev space $W^{1,p}$ with $1 \leq p < \infty$. By iterating the construction shown in Figure~\ref{f:disks}, in a self-similar manner, we obtain an example of the following:
\begin{itemize}
\item For any $1 \leq p < \infty$ there exist $u \in L^\infty([0,+\infty[ ; W^{1,p}(\T^d))$ and an initial configuration such that ${\rm mix}(\rho(t,\cdot)) \simeq e^{-Ct}$ for all $t\geq 0$.
\end{itemize}

The construction in~\cite{MR3904158} also reveals several properties of the regular Lagrangian flow associated with a Sobolev, but non-Lipschitz, velocity field. In particular, the flow of a velocity field in $L^\infty([0,T];W^{s,p}(\T^2))$, provided that $s$ and $p$ lie below the threshold for embedding into the class of Lipschitz velocity fields, can
\begin{itemize}
\item change the topology of sets (for instance, disconnect a connected set),
\item exhibit non-uniqueness of trajectories starting from a segment of initial points, and
\item compress such a segment to a single point (the regular Lagrangian flow has to preserve the two-dimensional Lebesgue measure ${\mathcal{L}}^2$, but we see here that it can compress the one-dimensional Hausdorff measure ${\mathcal{H}}^1$ to a delta measure; see also~\cite{MR4334730}). 
\end{itemize}	
A related construction was carried out in~\cite{MR3656475}. In that paper, the result is in some sense stronger: for any given initial configuration, one can construct a corresponding velocity field (depending on the initial data) that mixes it and belongs to $u \in L^\infty([0,+\infty[ ; W^{1,p}(\T^d))$. However, in certain regimes of large $p$, the decay of the mixing scale is slightly slower than exponential.

It is clear that the phenomena described above cannot occur for Lipschitz velocity fields. However, in~\cite{MR3904158}  it was shown that the exponential decay of the mixing scale can still be saturated by velocity fields that are uniformly-in-time Lipschitz in space. This was achieved through a variation of the self-similar scheme, adapted in~\cite{MR3904158} into a more general quasi–self-similar framework. The latter consists of replicating, at each integer time, a finite number of rescaled basic configurations, thereby avoiding the issues related to topological disconnection. This leads to the following result:
\begin{itemize}
\item There exist a smooth velocity field $u \in L^\infty([0,+\infty[ ; {\rm Lip}(\T^d))$ and an initial smooth configuration such that ${\rm mix}(\rho(t,\cdot)) \simeq e^{-Ct}$ for all $t\geq 0$.
\end{itemize}

This construction saturates the elementary bound obtained in Subsection~\ref{ss:lipPDE} and Section~\ref{s:flowestimates}, which follows from the Gr\"onwall inequality. It is easy to see, by means of a a linear velocity field, that Gr\"onwall’s estimate can be saturated pointwise, i.e., there exist pairs of trajectories whose distance increases exponentially. The result just discussed is stronger: it represents an integrated, or global, saturation of the Gr\"onwall inequality, in the sense that the exponential growth rate is attained at almost every point in space.

The preceding example further implies, via the interpolation estimate~\eqref{e:interpol}, that for smooth velocity fields and initial data, the positive Sobolev norms of the solution may increase exponentially in time. By suitably patching together such solutions and rescaling them to make the exponential growth arbitrarily fast, in~\cite{MR3933614} the following statement on loss of regularity for~\eqref{e:PDE} has been shown: 
\begin{itemize}
\item In dimension $d\geq 2$,  for any $1 \leq p < \infty$ there exist a bounded,  divergence-free velocity field $u \in L^\infty([0,T] ; W^{1,p}(\T^d))$ and a smooth initial datum $\bar\rho \in C^\infty(\T^d)$ such that for any $t>0$ and~$s>0$ the unique solution $\rho(t,\cdot) \not \in \dot{H}^s(\T^d)$. 
\end{itemize}
In fact, both the velocity field and the solution can be taken to be smooth for all times out of a fixed point in space. A variation of this result has been shown in~\cite{MR4430388}. This also implies loss of regularity for the regular Lagrangian flow of a non-Lipschitz velocity field, which has also been observed by means of a random construction in~\cite{MR3437603}. 

The above results, which stand in contrast to the DiPerna–Lions–Ambrosio theory~\cite{MR1022305,MR2096794}, suggest the following general principle.
For velocity fields in $W^{1,p}$, well-posedness of both~\eqref{e:PDE} and~\eqref{e:ODE} holds all the way down to $p = 1$, thereby extending the classical theory valid in the Lipschitz case $p = \infty$. However, the geometric and regularity properties of the flow and of the solution persist only at the level~$p = \infty$, that is, in the Lipschitz regime or nearby it. 

In the above works, the velocity field depends on the solution itself, and, due to the self-similar nature of the construction, it develops increasingly fine spatial scales as time progresses.
More recently, exploiting connections with the theory of dynamical systems, examples of universal mixers have been constructed in~\cite{MR4008523,MR4829551,MR4912975}. These are velocity fields that mix any initial configuration while maintaining uniform-in-time regularity in the spatial variable and periodicity in time.


\bibliographystyle{alpha}
\bibliography{mixref}

\begin{bibdiv}
\begin{biblist}

\bib{MR3268760}{article}{
      author={Alberti, Giovanni},
      author={Crippa, Gianluca},
      author={Mazzucato, Anna~L.},
       title={Exponential self-similar mixing and loss of regularity for
  continuity equations},
        date={2014},
        ISSN={1631-073X,1778-3569},
     journal={C. R. Math. Acad. Sci. Paris},
      volume={352},
      number={11},
       pages={901\ndash 906},
         url={https://doi.org/10.1016/j.crma.2014.08.021},
      review={\MR{3268760}},
}

\bib{MR3904158}{article}{
      author={Alberti, Giovanni},
      author={Crippa, Gianluca},
      author={Mazzucato, Anna~L.},
       title={Exponential self-similar mixing by incompressible flows},
        date={2019},
        ISSN={0894-0347,1088-6834},
     journal={J. Amer. Math. Soc.},
      volume={32},
      number={2},
       pages={445\ndash 490},
         url={https://doi.org/10.1090/jams/913},
      review={\MR{3904158}},
}

\bib{MR3933614}{article}{
      author={Alberti, Giovanni},
      author={Crippa, Gianluca},
      author={Mazzucato, Anna~L.},
       title={Loss of regularity for the continuity equation with
  non-{L}ipschitz velocity field},
        date={2019},
        ISSN={2524-5317,2199-2576},
     journal={Ann. PDE},
      volume={5},
      number={1},
       pages={Paper No. 9, 19},
         url={https://doi.org/10.1007/s40818-019-0066-3},
      review={\MR{3933614}},
}

\bib{MR2096794}{article}{
      author={Ambrosio, Luigi},
       title={Transport equation and {C}auchy problem for {$BV$} vector
  fields},
        date={2004},
        ISSN={0020-9910,1432-1297},
     journal={Invent. Math.},
      volume={158},
      number={2},
       pages={227\ndash 260},
         url={https://doi.org/10.1007/s00222-004-0367-2},
      review={\MR{2096794}},
}

\bib{MR3283066}{article}{
      author={Ambrosio, Luigi},
      author={Crippa, Gianluca},
       title={Continuity equations and {ODE} flows with non-smooth velocity},
        date={2014},
        ISSN={0308-2105,1473-7124},
     journal={Proc. Roy. Soc. Edinburgh Sect. A},
      volume={144},
      number={6},
       pages={1191\ndash 1244},
         url={https://doi.org/10.1017/S0308210513000085},
      review={\MR{3283066}},
}

\bib{MR4071413}{article}{
      author={Bianchini, Stefano},
      author={Bonicatto, Paolo},
       title={A uniqueness result for the decomposition of vector fields in
  {$\Bbb R^d$}},
        date={2020},
        ISSN={0020-9910,1432-1297},
     journal={Invent. Math.},
      volume={220},
      number={1},
       pages={255\ndash 393},
         url={https://doi.org/10.1007/s00222-019-00928-8},
      review={\MR{4071413}},
}

\bib{MR3862947}{article}{
      author={Bresch, Didier},
      author={Jabin, Pierre-Emmanuel},
       title={Global existence of weak solutions for compressible
  {N}avier-{S}tokes equations: thermodynamically unstable pressure and
  anisotropic viscous stress tensor},
        date={2018},
        ISSN={0003-486X,1939-8980},
     journal={Ann. of Math. (2)},
      volume={188},
      number={2},
       pages={577\ndash 684},
         url={https://doi.org/10.4007/annals.2018.188.2.4},
      review={\MR{3862947}},
}

\bib{MR2033003}{article}{
      author={Bressan, Alberto},
       title={An ill posed {C}auchy problem for a hyperbolic system in two
  space dimensions},
        date={2003},
        ISSN={0041-8994},
     journal={Rend. Sem. Mat. Univ. Padova},
      volume={110},
       pages={103\ndash 117},
      review={\MR{2033003}},
}

\bib{MR2033002}{article}{
      author={Bressan, Alberto},
       title={A lemma and a conjecture on the cost of rearrangements},
        date={2003},
        ISSN={0041-8994},
     journal={Rend. Sem. Mat. Univ. Padova},
      volume={110},
       pages={97\ndash 102},
      review={\MR{2033002}},
}

\bib{bressan2006prize}{misc}{
      author={Bressan, Alberto},
       title={Prize offered for the solution of a problem on mixing flows},
         how={\url{http://www.math.psu.edu/bressan/PSPDF/prize1.pdf}},
        date={2006},
        note={Accessed: 2025-09-30},
}

\bib{MR4377866}{article}{
      author={Bru\`e, Elia},
      author={Nguyen, Quoc-Hung},
       title={Sharp regularity estimates for solutions of the continuity
  equation drifted by {S}obolev vector fields},
        date={2021},
        ISSN={2157-5045,1948-206X},
     journal={Anal. PDE},
      volume={14},
      number={8},
       pages={2539\ndash 2559},
         url={https://doi.org/10.2140/apde.2021.14.2539},
      review={\MR{4377866}},
}

\bib{MR4263701}{article}{
      author={Bru\`e, Elia},
      author={Nguyen, Quoc-Hung},
       title={Sobolev estimates for solutions of the transport equation and
  {ODE} flows associated to non-{L}ipschitz drifts},
        date={2021},
        ISSN={0025-5831,1432-1807},
     journal={Math. Ann.},
      volume={380},
      number={1-2},
       pages={855\ndash 883},
         url={https://doi.org/10.1007/s00208-020-01988-5},
      review={\MR{4263701}},
}

\bib{MR4797689}{article}{
      author={Bru\`e, Elia},
      author={Zelati, Michele~Coti},
      author={Marconi, Elio},
       title={Enhanced dissipation for two-dimensional {H}amiltonian flows},
        date={2024},
        ISSN={0003-9527,1432-0673},
     journal={Arch. Ration. Mech. Anal.},
      volume={248},
      number={5},
       pages={Paper No. 84, 37},
         url={https://doi.org/10.1007/s00205-024-02034-3},
      review={\MR{4797689}},
}

\bib{MR4733112}{article}{
      author={Coti~Zelati, Michele},
      author={Crippa, Gianluca},
      author={Iyer, Gautam},
      author={Mazzucato, Anna~L.},
       title={Mixing in incompressible flows: transport, dissipation, and their
  interplay},
        date={2024},
        ISSN={0002-9920,1088-9477},
     journal={Notices Amer. Math. Soc.},
      volume={71},
      number={5},
       pages={593\ndash 604},
      review={\MR{4733112}},
}

\bib{MR2369485}{article}{
      author={Crippa, Gianluca},
      author={De~Lellis, Camillo},
       title={Estimates and regularity results for the {D}i{P}erna-{L}ions
  flow},
        date={2008},
        ISSN={0075-4102,1435-5345},
     journal={J. Reine Angew. Math.},
      volume={616},
       pages={15\ndash 46},
         url={https://doi.org/10.1515/CRELLE.2008.016},
      review={\MR{2369485}},
}

\bib{MR4430388}{article}{
      author={Crippa, Gianluca},
      author={Elgindi, Tarek},
      author={Iyer, Gautam},
      author={Mazzucato, Anna~L.},
       title={Growth of {S}obolev norms and loss of regularity in transport
  equations},
        date={2022},
        ISSN={1364-503X,1471-2962},
     journal={Philos. Trans. Roy. Soc. A},
      volume={380},
      number={2225},
       pages={Paper No. 24, 12},
         url={https://doi.org/10.1098/rsta.2021.0024},
      review={\MR{4430388}},
}

\bib{MR3703559}{article}{
      author={Crippa, Gianluca},
      author={Schulze, Christian},
       title={Cellular mixing with bounded palenstrophy},
        date={2017},
        ISSN={0218-2025,1793-6314},
     journal={Math. Models Methods Appl. Sci.},
      volume={27},
      number={12},
       pages={2297\ndash 2320},
         url={https://doi.org/10.1142/S0218202517500452},
      review={\MR{3703559}},
}

\bib{MR2009116}{article}{
      author={Depauw, Nicolas},
       title={Non unicit\'e{} des solutions born\'ees pour un champ de vecteurs
  {BV} en dehors d'un hyperplan},
        date={2003},
        ISSN={1631-073X,1778-3569},
     journal={C. R. Math. Acad. Sci. Paris},
      volume={337},
      number={4},
       pages={249\ndash 252},
         url={https://doi.org/10.1016/S1631-073X(03)00330-3},
      review={\MR{2009116}},
}

\bib{MR1022305}{article}{
      author={DiPerna, Ronald~J.},
      author={Lions, Pierre-Louis},
       title={Ordinary differential equations, transport theory and {S}obolev
  spaces},
        date={1989},
        ISSN={0020-9910,1432-1297},
     journal={Invent. Math.},
      volume={98},
      number={3},
       pages={511\ndash 547},
         url={https://doi.org/10.1007/BF01393835},
      review={\MR{1022305}},
}

\bib{MR3824693}{incollection}{
      author={Doering, Charles~R.},
      author={Nobili, Camilla},
       title={Lectures on stirring, mixing and transport},
        date={2017},
   booktitle={Transport, fluids, and mixing},
      series={Partial Differ. Equ. Meas. Theory},
   publisher={De Gruyter Open, Warsaw},
       pages={8\ndash 34},
         url={https://doi.org/10.1515/9783110571240-005},
      review={\MR{3824693}},
}

\bib{MR4829551}{article}{
      author={Elgindi, Tarek~M.},
      author={Liss, Kyle},
       title={Norm growth, non-uniqueness, and anomalous dissipation in passive
  scalars},
        date={2024},
        ISSN={0003-9527,1432-0673},
     journal={Arch. Ration. Mech. Anal.},
      volume={248},
      number={6},
       pages={Paper No. 120, 28},
         url={https://doi.org/10.1007/s00205-024-02056-x},
      review={\MR{4829551}},
}

\bib{MR4912975}{article}{
      author={Elgindi, Tarek~M.},
      author={Liss, Kyle},
      author={Mattingly, Jonathan~C.},
       title={Optimal enhanced dissipation and mixing for a time-periodic,
  {L}ipschitz velocity field on {$\Bbb{T}^2$}},
        date={2025},
        ISSN={0012-7094,1547-7398},
     journal={Duke Math. J.},
      volume={174},
      number={7},
       pages={1209\ndash 1260},
         url={https://doi.org/10.1215/00127094-2024-0057},
      review={\MR{4912975}},
}

\bib{MR4008523}{article}{
      author={Elgindi, Tarek~M.},
      author={Zlato\v{s}, Andrej},
       title={Universal mixers in all dimensions},
        date={2019},
        ISSN={0001-8708,1090-2082},
     journal={Adv. Math.},
      volume={356},
       pages={106807, 33},
         url={https://doi.org/10.1016/j.aim.2019.106807},
      review={\MR{4008523}},
}

\bib{MR4334730}{article}{
      author={Fefferman, Charles~L.},
      author={Pooley, Benjamin~C.},
      author={Rodrigo, Jos\'e~L.},
       title={Non-conservation of dimension in divergence-free solutions of
  passive and active scalar systems},
        date={2021},
        ISSN={0003-9527,1432-0673},
     journal={Arch. Ration. Mech. Anal.},
      volume={242},
      number={3},
       pages={1445\ndash 1478},
         url={https://doi.org/10.1007/s00205-021-01708-6},
      review={\MR{4334730}},
}

\bib{2402.11642}{misc}{
      author={Huysmans, Lucas},
      author={Said, Ayman~Rimah},
       title={Quantitative estimates in passive scalar transport: A unified
  approach in {$W^{1,p}$} via {C}hrist-{J}ourné commutator estimates},
        date={2024},
}

\bib{2504.03023}{misc}{
      author={Huysmans, Lucas},
      author={Said, Ayman~Rimah},
       title={Mixing estimates for passive scalar transport by {$BV$} vector
  fields},
        date={2025},
}

\bib{MR3207161}{article}{
      author={Iyer, Gautam},
      author={Kiselev, Alexander},
      author={Xu, Xiaoqian},
       title={Lower bounds on the mix norm of passive scalars advected by
  incompressible enstrophy-constrained flows},
        date={2014},
        ISSN={0951-7715,1361-6544},
     journal={Nonlinearity},
      volume={27},
      number={5},
       pages={973\ndash 985},
         url={https://doi.org/10.1088/0951-7715/27/5/973},
      review={\MR{3207161}},
}

\bib{MR3437603}{article}{
      author={Jabin, Pierre-Emmanuel},
       title={Critical non-{S}obolev regularity for continuity equations with
  rough velocity fields},
        date={2016},
        ISSN={0022-0396,1090-2732},
     journal={J. Differential Equations},
      volume={260},
      number={5},
       pages={4739\ndash 4757},
         url={https://doi.org/10.1016/j.jde.2015.11.028},
      review={\MR{3437603}},
}

\bib{MR3799259}{article}{
      author={L\'eger, Flavien},
       title={A new approach to bounds on mixing},
        date={2018},
        ISSN={0218-2025,1793-6314},
     journal={Math. Models Methods Appl. Sci.},
      volume={28},
      number={5},
       pages={829\ndash 849},
         url={https://doi.org/10.1142/S0218202518500215},
      review={\MR{3799259}},
}

\bib{MR2801050}{article}{
      author={Lin, Zhi},
      author={Thiffeault, Jean-Luc},
      author={Doering, Charles~R.},
       title={Optimal stirring strategies for passive scalar mixing},
        date={2011},
        ISSN={0022-1120,1469-7645},
     journal={J. Fluid Mech.},
      volume={675},
       pages={465\ndash 476},
         url={https://doi.org/10.1017/S0022112011000292},
      review={\MR{2801050}},
}

\bib{MR3026556}{article}{
      author={Lunasin, Evelyn},
      author={Lin, Zhi},
      author={Novikov, Alexei},
      author={Mazzucato, Anna},
      author={Doering, Charles~R.},
       title={Optimal mixing and optimal stirring for fixed energy, fixed
  power, or fixed palenstrophy flows},
        date={2012},
        ISSN={0022-2488,1089-7658},
     journal={J. Math. Phys.},
      volume={53},
      number={11},
       pages={115611, 15},
         url={https://doi.org/10.1063/1.4752098},
      review={\MR{3026556}},
}

\bib{MR2200439}{article}{
      author={Mathew, George},
      author={Mezi\'c, Igor},
      author={Petzold, Linda},
       title={A multiscale measure for mixing},
        date={2005},
        ISSN={0167-2789,1872-8022},
     journal={Phys. D},
      volume={211},
      number={1-2},
       pages={23\ndash 46},
         url={https://doi.org/10.1016/j.physd.2005.07.017},
      review={\MR{2200439}},
}

\bib{MR4747245}{article}{
      author={Meyer, David},
      author={Seis, Christian},
       title={Propagation of regularity for transport equations. {A}
  {L}ittlewood-{P}aley approach},
        date={2024},
        ISSN={0022-2518,1943-5258},
     journal={Indiana Univ. Math. J.},
      volume={73},
      number={2},
       pages={445\ndash 472},
      review={\MR{4747245}},
}

\bib{MR4242824}{article}{
      author={Nguyen, Quoc-Hung},
       title={Quantitative estimates for regular {L}agrangian flows with {$BV$}
  vector fields},
        date={2021},
        ISSN={0010-3640,1097-0312},
     journal={Comm. Pure Appl. Math.},
      volume={74},
      number={6},
       pages={1129\ndash 1192},
         url={https://doi.org/10.1002/cpa.21992},
      review={\MR{4242824}},
}

\bib{MR3141856}{article}{
      author={Seis, Christian},
       title={Maximal mixing by incompressible fluid flows},
        date={2013},
        ISSN={0951-7715,1361-6544},
     journal={Nonlinearity},
      volume={26},
      number={12},
       pages={3279\ndash 3289},
         url={https://doi.org/10.1088/0951-7715/26/12/3279},
      review={\MR{3141856}},
}

\bib{MR2876867}{article}{
      author={Thiffeault, Jean-Luc},
       title={Using multiscale norms to quantify mixing and transport},
        date={2012},
        ISSN={0951-7715,1361-6544},
     journal={Nonlinearity},
      volume={25},
      number={2},
       pages={R1\ndash R44},
         url={https://doi.org/10.1088/0951-7715/25/2/R1},
      review={\MR{2876867}},
}

\bib{MR3656475}{article}{
      author={Yao, Yao},
      author={Zlato\v{s}, Andrej},
       title={Mixing and un-mixing by incompressible flows},
        date={2017},
        ISSN={1435-9855,1435-9863},
     journal={J. Eur. Math. Soc. (JEMS)},
      volume={19},
      number={7},
       pages={1911\ndash 1948},
         url={https://doi.org/10.4171/JEMS/709},
      review={\MR{3656475}},
}

\bib{MR4026549}{article}{
      author={Zillinger, Christian},
       title={On geometric and analytic mixing scales: comparability and
  convergence rates for transport problems},
        date={2019},
        ISSN={2578-5885,2578-5893},
     journal={Pure Appl. Anal.},
      volume={1},
      number={4},
       pages={543\ndash 570},
         url={https://doi.org/10.2140/paa.2019.1.543},
      review={\MR{4026549}},
}

\end{biblist}
\end{bibdiv}

\end{document}